\documentclass[11pt]{article}
\usepackage{amssymb,amsmath,latexsym}
\usepackage{pdfsync}

\allowdisplaybreaks

\begin{document}

\begin{doublespace}

\def\1{{\bf 1}}
\def\ind{{\bf 1}}
\def\nn{\nonumber}

\def\sA {{\cal A}} \def\sB {{\cal B}} \def\sC {{\cal C}}
\def\sD {{\cal D}} \def\sE {{\cal E}} \def\sF {{\cal F}}
\def\sG {{\cal G}} \def\sH {{\cal H}} \def\sI {{\cal I}}
\def\sJ {{\cal J}} \def\sK {{\cal K}} \def\sL {{\cal L}}
\def\sM {{\cal M}} \def\sN {{\cal N}} \def\sO {{\cal O}}
\def\sP {{\cal P}} \def\sQ {{\cal Q}} \def\sR {{\cal R}}
\def\sS {{\cal S}} \def\sT {{\cal T}} \def\sU {{\cal U}}
\def\sV {{\cal V}} \def\sW {{\cal W}} \def\sX {{\cal X}}
\def\sY {{\cal Y}} \def\sZ {{\cal Z}}

\def\bA {{\mathbb A}} \def\bB {{\mathbb B}} \def\bC {{\mathbb C}}
\def\bD {{\mathbb D}} \def\bE {{\mathbb E}} \def\bF {{\mathbb F}}
\def\bG {{\mathbb G}} \def\bH {{\mathbb H}} \def\bI {{\mathbb I}}
\def\bJ {{\mathbb J}} \def\bK {{\mathbb K}} \def\bL {{\mathbb L}}
\def\bM {{\mathbb M}} \def\bN {{\mathbb N}} \def\bO {{\mathbb O}}
\def\bP {{\mathbb P}} \def\bQ {{\mathbb Q}} \def\bR {{\mathbb R}}
\def\bS {{\mathbb S}} \def\bT {{\mathbb T}} \def\bU {{\mathbb U}}
\def\bV {{\mathbb V}} \def\bW {{\mathbb W}} \def\bX {{\mathbb X}}
\def\bY {{\mathbb Y}} \def\bZ {{\mathbb Z}}
\def\R {{\mathbb R}} \def\RR {{\mathbb R}} \def\H {{\mathbb H}}
\def\n{{\bf n}} \def\Z {{\mathbb Z}}

\newcommand{\expr}[1]{\left( #1 \right)}
\newcommand{\cl}[1]{\overline{#1}}
\newtheorem{thm}{Theorem}[section]
\newtheorem{lemma}[thm]{Lemma}
\newtheorem{defn}[thm]{Definition}
\newtheorem{prop}[thm]{Proposition}
\newtheorem{corollary}[thm]{Corollary}
\newtheorem{remark}[thm]{Remark}
\newtheorem{example}[thm]{Example}
\numberwithin{equation}{section}
\def\ee{\varepsilon}
\def\qed{{\hfill $\Box$ \bigskip}}
\def\NN{{\cal N}}
\def\AA{{\cal A}}
\def\MM{{\cal M}}
\def\BB{{\cal B}}
\def\CC{{\cal C}}
\def\LL{{\cal L}}
\def\DD{{\cal D}}
\def\FF{{\cal F}}
\def\EE{{\cal E}}
\def\QQ{{\cal Q}}
\def\RR{{\mathbb R}}
\def\R{{\mathbb R}}
\def\L{{\bf L}}
\def\K{{\bf K}}
\def\S{{\bf S}}
\def\A{{\bf A}}
\def\E{{\mathbb E}}
\def\F{{\bf F}}
\def\P{{\mathbb P}}
\def\N{{\mathbb N}}
\def\eps{\varepsilon}
\def\wh{\widehat}
\def\wt{\widetilde}
\def\pf{\noindent{\bf Proof.} }
\def\pff{\noindent{\bf Proof} }
\def\cp{\mathrm{Cap}}

\title{\Large \bf Martin boundary for some symmetric L\'evy processes}

\author{{\bf Panki Kim}\thanks{
This work was supported by the National Research Foundation of Korea(NRF) grant funded by the Korea government(MEST) (NRF-2013R1A2A2A01004822)
}
\quad {\bf Renming Song\thanks{Research supported in part by a grant from the Simons Foundation (208236)}} \quad and
\quad {\bf Zoran Vondra\v{c}ek}
}

\date{}

\maketitle

\begin{abstract}
In this paper we study the Martin boundary of open sets with respect to a large class of
purely discontinuous symmetric L\'evy processes in $\R^d$. We show that, if
$D\subset \R^d$ is an open set which is $\kappa$-fat at a boundary point $Q\in \partial D$,
then there is exactly one Martin boundary point associated with $Q$ and this
Martin boundary point is minimal.
\end{abstract}



\section{Introduction}\label{sec:1}

The Martin boundary of an open set $D$ is an abstract boundary introduced in 1941
by Martin \cite{martin} so that every nonnegative
classical harmonic function in $D$ can be written as an integral of the Martin kernel with respect to a finite
measure on the Martin boundary. This integral representation is called a Martin representation.
The concepts of Martin boundary and Martin kernel were extended to general
Markov processes by Kunita and Watanabe \cite{KW} in 1965.
In order for the Martin representation to be useful, one needs to
have a better understanding of the Martin boundary, for instance, its relation with the Euclidean boundary.
In 1970, Hunt and Wheeden \cite{HW} proved that, in the classical case, the Martin boundary of
a bounded Lipschitz domain coincides with its Euclidean boundary.
Subsequently, a lot of progress
has been made in studying the Martin boundary in the classical case.

With the help of the boundary Harnack principle for rotationally invariant $\alpha$-stable ($\alpha\in (0, 2)$)
processes established in \cite{B}, it was proved in \cite{B99, CS98, MS00} that the Martin boundary,
with respect to the rotationally invariant $\alpha$-stable process, of a bounded Lipschitz domain $D$
coincides with its Euclidean boundary and that any nonnegative harmonic function with respect to the killed
rotationally invariant $\alpha$-stable process in $D$ can be written uniquely as an integral of the Martin
kernel with respect to a finite measure on $\partial D$.
In \cite{SW99} this result was extended
to bounded $\kappa$-fat open sets. The Martin boundary, with respect to
truncated stable processes, of any roughly connected $\kappa$-fat open set was shown in \cite{KS08} to coincide with
its Euclidean boundary. In \cite{KSV1}, the results of \cite{SW99} were extended to a large class of purely
discontinuous subordinate Brownian motions. For Martin boundary at infinity with respect to subordinate Brownian motions,
see \cite{KSV9}.

In this paper, we study the Martin boundary of open set $D\subset \R^d$ with respect to a large class of symmetric, not
necessarily rotationally invariant, (transient) L\'evy processes killed upon exiting $D$.
We show that if $D$ is an open set and $D$ is $\kappa$-fat at a single point $Q\in \partial D$,
then the Martin boundary associated with $Q$ consists of exactly one point and the corresponding
Martin kernel is a minimal harmonic function.
Another point is that, unlike \cite{B99, CS98, KS08, KSV1, MS00, SW99}, the set $D$ is not necessarily bounded.
In the case when $D$ is unbounded,
we do not study the Martin boundary associated with infinite boundary points.

Now we describe the class of processes we are going to work with.

Throughout this paper,
$r\mapsto j(r)$ is a strictly positive and non-increasing function
on $(0, \infty)$ satisfying
\begin{equation}\label{e:j-decay}
j(r)\le cj(r+1) \qquad \hbox{for } r\ge 1,
\end{equation}
and  $X=(X_t,\P_x)$ is a purely discontinuous symmetric L\'evy process with L\'evy exponent $\Psi_X(\xi)$ so that
$$
\E_x\left[e^{i\xi\cdot(X_t-X_0)}\right]=e^{-t\Psi_X(\xi)}, \quad \quad t>0, x\in \R^d, \xi\in \R^d.
$$
We assume that the L\'evy measure of $X$ has a density $J_X$ such that
\begin{equation}\label{e:psi1}
\gamma^{-1}_1 j(|y|)\le J_X(y) \le \gamma_1 j(|y|), \quad \mbox{for all } y\in \R^d\, ,
\end{equation}
for some $\gamma_1 >1$.
Since
$
\int_0^\infty j(r) (1\wedge r^2) r^{d-1}dr < \infty$ by \eqref{e:psi1},
the function $x \to j(|x|)$ is the L\'evy density of  an isotropic unimodal L\'evy process whose characteristic exponent is
\begin{equation}\label{e:Psi}
\Psi(|\xi|)= \int_{\R^d}(1-\cos(\xi\cdot y))j(|y|)dy.
\end{equation}
The L\'evy exponent $\Psi_X$ can be written as
$$
\Psi_X(\xi)= \int_{\R^d}(1-\cos(\xi\cdot y))J_X(y)dy
$$
and, clearly by \eqref{e:psi1}, it satisfies
\begin{equation}\label{e:psi21}
\gamma^{-1}_1 \Psi(|\xi|)\le \Psi_X(\xi) \le \gamma_1 \Psi(|\xi|),
\quad \mbox{for all } \xi\in \R^d\, .
\end{equation}
 The function $\Psi$ may be not increasing.
However, if we put $\Psi^*(r):= \sup_{s \le r} \Psi(s)$, then,
by \cite[Proposition 2]{BGR1}  (cf.~also \cite[Proposition 1]{G}), we have
$$
\Psi(r) \le\Psi^*(r) \le \pi^2 \Psi(r).
$$
Thus by \eqref{e:psi21},
\begin{equation}\label{e:psi3}
(\pi^2\gamma_1)^{-1}  \Psi^*(|\xi|)\le \Psi_X(\xi) \le \gamma_1 \Psi^*(|\xi|),
\quad \mbox{for all } \xi  \in \R^d\, .
\end{equation}

We will always assume that $\Psi$ satisfies the following scaling condition at infinity:

\medskip
\noindent
({\bf H}):
There exist constants $0<\delta_1\le \delta_2 <1$ and $a_1, a_2>0$  such that
\begin{equation}\label{e:H1n}
a_1\lambda^{2\delta_1} \Psi(t) \le \Psi(\lambda t) \le a_2 \lambda^{2\delta_2} \Psi(t), \quad \lambda \ge 1, t \ge 1\, .
\end{equation}
Then by  \cite[(15) and Corollary 22]{BGR1},
for every $R>0$, there exists $c=c(R)>1$ such that
\begin{equation}\label{e:asmpbofjat0n}
c^{-1}\frac{\Psi(r^{-1})}{r^d} \le j(r)\le c \frac{\Psi(r^{-1})}{r^d}
\quad \hbox{for } r\in (0, R].
\end{equation}
Note that the class of purely discontinuous symmetric L\'evy processes considered in this paper contains some of
the purely discontinuous isotropic unimodal L\'evy processes dealt with in \cite{BGR1}.

Let us now formulate precisely the main result of this paper.
\begin{defn}\label{d:kappa-fat-at-point}
Let $D\subset \R^d$ be an open set and $Q\in \R^d$. We say that $D$ is $\kappa$-fat
at $Q$ for some $\kappa\in (0, \frac12)$, if there exists $R>0$ such that for all $r\in
(0,R]$, there is a ball  $B(A_r(Q) , \kappa r)  \subset D\cap B(Q,r)$.
The pair $(R, \kappa)$ is called the characteristics of
the $\kappa$-fat open set $D$ at $Q$. We say that an open set $D$
is $\kappa$-fat with characteristics $(R, \kappa)$ if $D$ is $\kappa$-fat
at $Q \in \partial D$ with characteristics $(R, \kappa)$ for all $Q \in  \partial D$ .
\end{defn}

For $D\subset \R^d$ we denote by $\partial_M D$ the Martin boundary of $D$. A point $w\in \partial_MD$ is said to be associated
with $Q$ if there is a sequence $(y_n)_{n\ge 1}\subset D$ converging to $w$ in the Martin topology and to $Q$ in the Euclidean topology. The set of Martin boundary points associated with $Q$ is denoted by $\partial_M^QD$.

\begin{thm}\label{t:main-theorem}
Suppose that the assumption {\bf (H)} is satisfied.
Let $X$ be a symmetric L\'evy process
with a  L\'evy density satisfying  \eqref{e:psi1}
and let $D$ be an open subset of $\R^d$
which is $\kappa$-fat at $Q\in \partial D$. If $D$ is bounded, then $\partial_M^QD$
consists of exactly one point and this point is a minimal Martin boundary point.
If $D$ is unbounded and the process $X$ is transient, the same conclusion is true.
\end{thm}

In the case when $D$ is unbounded,
a natural assumption would be that $D$ is Greenian, that is,
the killed process $X^D$ is transient. Unfortunately, under the assumption {\bf (H)},
which governs the behavior of the process in small space, it seems difficult to obtain
estimates of the Green function $G_D(x,y)$ when either $x$ or $y$ becomes large. This
is why in the case of unbounded $D$ we impose the transience assumption on $X$ which gives
the asymptotic behavior of the free Green function $G(x,y)$, cf.~Lemma \ref{l:gfnatinfty}.

We begin the paper by showing that only small modifications are needed to extend some
results from the isotropic case studied in \cite{KSV2}
to the symmetric L\'evy processes $X$ considered in this paper.
These results include exit time estimates, Poisson kernel estimates and Harnack inequality.
A little more work is needed to establish the upper and the lower bound on the Green
function $G_D$. Those are used to obtain sharp bounds on the Poisson kernel and the
boundary Harnack principle in the same way as in \cite{KSV2}. In Section 3
we follow the well-established route, see \cite{B, KSV1, KSV6, KSV9}, to identify
the Martin boundary point associated to $Q$. After preliminary estimates about
harmonic functions, we first show that the oscillation reduction lemma,
see \cite[Lemma 16]{B}, is valid in our setting (with essentially the same proof).
The lemma almost immediately implies that $\partial_M^Q D$ consists of
exactly one point.
We then show that this point is a minimal Martin boundary point.
We end the paper by giving the Martin representation for bounded $\kappa$-fat
open sets.

We finish this introduction by setting up some notation and conventions.
We use ``$:=$" to denote a definition, which
is  read as ``is defined to be"; we denote $a \wedge b := \min \{ a, b\}$,
$a \vee b := \max \{ a, b\}$;
we denote by $B(x, r)$ the open
ball centered at $x\in \bR^d$ with radius $r>0$;
for any two positive functions $f$ and $g$,
$f\asymp g$ means that there is a positive constant $c\geq 1$
so that $c^{-1}\, g \leq f \leq c\, g$ on their common domain of
definition;
for any Borel subset $E\subset\bR^d$ and $x\in E$,
${\rm diam}(E)$ stands for the diameter of $E$ and
$\delta_E(x)$ stands for the Euclidean distance between
$x$ and $E^c$; $\N$ is the set of natural numbers.

In this paper, we use the following convention: The values of
the constants  $R, \delta_1, \delta_2, C_1, C_2, C_3, C_4$ remain the same
throughout this paper, while $c, c_0, c_1, c_2, \ldots$  represent constants
whose values are unimportant and may change. All constants are positive finite numbers.
The labeling of the constants $c_0, c_1, c_2, \ldots$ starts anew in the statement and proof of each result.
The dependence of constant $c$ on dimension $d$ is not mentioned explicitly.

\section{Green function estimates}

Let $S=(S_t:\, t\ge 0)$ be a subordinator with no drift.
The Laplace exponent $\phi$ of $S$ is a Bernstein function and admits the following
representation
$$
\phi(\lambda)=\int^\infty_0(1-e^{-\lambda t})\mu(dt), \qquad \lambda>0,
$$
where $\mu$ satisfies $\int^\infty_0(1\wedge t)\mu(dt)<\infty$. $\mu$ is called the
L\'evy measure of $S$ or $\phi$.
The function $\phi$ is called a complete Bernstein function if $\mu$ has a completely monotone density.

The following elementary result observed in \cite{KSV6, KSV8} will be used several times later in this paper.

\begin{lemma}\label{l:phi-property}
If $\phi$ is a Bernstein function, then for all  $\lambda, t>0$,
$1 \wedge  \lambda\le {\phi(\lambda t)}/{\phi(t)} \le 1 \vee
 \lambda$.
\end{lemma}

Recall that,
throughout this paper, we assume that
$X=(X_t, \P_x)$ is a purely discontinuous symmetric
L\'evy process in $\R^d$ with L\'evy exponent $\Psi_X(\xi)$ and
a L\'evy density $J_X$ satisfying \eqref{e:psi1}.

It follows from
\cite[(28)]{BGR1}, \eqref{e:psi21},
\eqref{e:H1n} and \eqref{e:asmpbofjat0n}
that there exist a constant $\gamma_2>1$ and a complete Bernstein function $\phi$
such that
\begin{equation}\label{e:psi2}
\gamma_2^{-1} \phi(|\xi|^2)\le \Psi(|\xi|) \le \gamma_2 \phi(|\xi|^2),
\quad \mbox{for all }
\xi \in \R^d\, ,
\end{equation}
and $j$
enjoys the following property:
for every $R>0$,
\begin{equation}\label{e:asmpbofjat0}
j(r)\asymp \frac{\phi(r^{-2})}{r^d}
\quad \hbox{for } r\in (0, R].
\end{equation}
Furthermore,
there exist  $b_1, b_2>0$  such that
\begin{equation}\label{e:H1}
b_1\lambda^{\delta_1} \phi(t) \le \phi(\lambda t) \le b_2 \lambda^{\delta_2} \phi(t), \quad \lambda \ge 1, t \ge 1\, .
\end{equation}
Throughout this paper, we assume that
$\phi$ is the above complete Bernstein function.

From Lemma \ref{l:phi-property} and \eqref{e:asmpbofjat0}, we also get that
for every $R>0$,
\begin{equation}\label{e:doubling-condition}
j(r)\le c j(2r)\, ,\quad r\in (0, R]\, .
\end{equation}

The infinitesimal generator $\L$ of $X$
is given by
\begin{equation}\label{3.1}
\L f(x)=\int_{\R^d}\left( f(x+y)-f(x)-y\cdot \nabla f(x)
{\bf 1}_{\{|y|\le1\}}
 \right)\, J_X(y)dy
\end{equation}
for $f\in C_b^2(\R^d)$. Furthermore, for every $f\in C_b^2(\R^d)$,
$
f(X_t)-f(X_0)-\int_0^t \L f(X_s)\, ds
$
is a $\P_x$-martingale for every $x\in \R^d$.

The following two results are valid without assuming ({\bf H}).
The next lemma is a special case of \cite[Corollary 1]{G}.

\begin{lemma}\label{l:j-upper}
 There exists a constant $c>0$ depending only on $d$ such that
    \begin{equation}\label{e:j-up}
c^{-1} \Psi(r^{-1}) \le  \frac{1}{r^2} \int_0^r s^{d+1} j(s)ds + \int_r^\infty  s^{d-1} j(s)ds\le  c \Psi(r^{-1})\, ,\qquad \forall r > 0.
    \end{equation}
\end{lemma}

\begin{lemma}\label{l:lnew}
There exists a constant $c=c(\Psi, \gamma_1, \gamma_2)>0$ such that for every
$f\in C^2_b(\R^d)$ with $0\leq f \leq 1$,
$$
\L f_r(x) \le   c\,
 \phi(r^{-2})
 \left( 2+ \sup_{y}\sum_{j,k} |(\partial^2/\partial y_j\partial y_k) f(y)|
 \right),
 \quad \text{for every } x \in \R^d,
r>0,
$$
where $f_r(y):=f(y/r)$.
\end{lemma}

\pf
Using \eqref{e:psi1} and  Lemma \ref{l:j-upper} or \cite[Corollary 3]{BGR1},
 this result can be obtained
by following the proof of \cite[Lemma 4.2]{KSV7}. We omit the details.
\qed

For any open set $D$, we use $\tau_D$ to denote the first exit
time of $D$, i.e., $\tau_D=\inf\{t>0: \, X_t\notin D\}$.

Using Lemmas \ref{l:phi-property} and \ref{l:lnew} and  \eqref{e:psi2},
the proof of the next result is the same as those of
\cite[Lemmas 13.4.1 and 13.4.2]{KSV3}. Thus we omit the proof.

\begin{lemma}\label{L3.2}
There exists a constant $c=c(\Psi,\gamma_1, \gamma_2)>0$ such that
for every
$r>0$ and every $x\in \R^d$,
$$
\inf_{z\in B(x,r/2)} \E_z \left[\tau_{B(x,r)} \right] \geq
\frac{c}{\phi(r^{-2})}.
$$
\end{lemma}

The idea of the following key result comes from \cite{Sz}.
\begin{lemma}\label{l:tau}
There exists a constant $c=c(\Psi, \gamma_1, \gamma_2)>0$ such that for any
$r>0$ and $x_0 \in \R^d$,
\begin{eqnarray*}
\E_x[\tau_{B(x_0,r)}]\le  c\, (\phi(r^{-2})\phi((r-|x-x_0|)^{-2}))^{-1/2}\, ,
\qquad x\in B(x_0, r).
\end{eqnarray*}
\end{lemma}

\pf Without loss of generality, we may assume that $x_0=0$. We fix
$x\neq 0$ and  put $Z_t=\frac{X_t\cdot x}{|x|}$.
Then,
$Z_t$ is a one dimensional symmetric L\'evy process in $\R$ with L\'evy exponent
$$
\Psi_Z(\theta)=\int_{\R^d}\left(1-\cos(\frac{\theta x}{|x|}\cdot y)\right)J_X(y)dy,
\qquad \theta\in \R.
$$
By \eqref{e:psi1},
$$
\Psi_Z(\theta) \asymp \int_{\R^d}\left(1-\cos(\frac{\theta x}{|x|}\cdot y)\right)j(|y|)dy=\Psi(\theta).
$$

It is easy to see that, if
$X_t\in B(0, r)$, then $|Z_t|<r$, hence
$\E_x[\tau_{B(0, r)}]\le \E_{|x|}[\tilde \tau],$
where $\tilde \tau=\inf\{t>0: |Z_t|\ge r\}$.
By \cite[(2.17)]{BGR2},  the proof of \cite[Proposition 2.4]{BGR2} and Lemma \ref{l:j-upper}.
$$\E_x[\tau_{B(0, r)}]\le \E_{|x|}[\tilde \tau] \le c\, (\Psi(r^{-1})
\Psi((r-|x|)^{-1}))^{-1/2}.
$$
Now the assertion of the lemma follows immediately by \eqref{e:psi2}.
\qed

Given  an open set $D\subset \R^d$, we define $X^D_t(\omega)=X_t(\omega)$ if $t< \tau_D(\omega)$ and
$X^D_t(\omega)=\partial$ if $t\geq  \tau_D(\omega)$, where $\partial$ is a cemetery state.
We now recall the definitions of harmonic functions with respect to $X$ and with respect to $X^D$.

\begin{defn}\label{def:har1}
Let $D$ be an open subset of $\R^d$.
A nonnegative function $u$ on $\R^d$ is said to be
\begin{description}
\item{\rm{(1)}}  harmonic in $D$ with respect to $X$ if
$$
u(x)= \E_x\left[u(X_{\tau_{
U}})\right],
\qquad x\in
U,
$$
for every open set $
U$ whose closure is a compact
subset of $D$;

\item{\rm{(2)}}
regular harmonic in $D$ with respect to $X$ if
for each $x \in D$,
$$
u(x)= \E_x\left[u(X_{\tau_{D}}); \tau_D<\infty\right].
$$
\end{description}
\end{defn}

\begin{defn}\label{def:harkilled}
Let $D$ be an open subset of $\R^d$.
A nonnegative function $u$ on $D$ is said to be
harmonic with respect to $X^D$ if
$$
u(x)= \E_x\left[u(X_{\tau_{
U}})\right],
\qquad x\in
U,
$$
for every open set $
U$ whose closure is a compact
subset of $D$.
\end{defn}

Obviously, if $u$ is harmonic with respect to $X^D$, then the function which is equal to
$u$ in $D$ and zero outside $D$ is harmonic with respect to $X$ in $D$.

Since our $X$ satisfies \cite[(1.6), {\bf (UJS)}]{CKK2}, by \cite[Theorem 1.4]{CKK2}
and  using the standard chain argument one has  the following form of Harnack inequality.

\begin{thm}\label{hi}
For every $a \in (0,1)$, there exists $c=c(a, \Psi, \gamma_1, \gamma_2)>0$ such that
for every $r \in (0, 1)$,
$x_0 \in {\mathbb R}^d$, and every function $u$
which is nonnegative on ${\mathbb R}^d$
and harmonic with respect to $X$ in $B(x_0, r)$, we have
$$
u(x)\,\le \,c\, u(y), \quad \textrm{for all }x, y\in B(x_0, ar)\, .
$$
\end{thm}

Let $D\subset \R^d$ be an open set.
Since $J_X$
satisfies the assumption \cite[(1.6)]{CKK2}, by \cite[Theorem 3.1]{CKK2}, bounded functions that
are harmonic in $D$ with respect to $X$ are H\"older continuous.
Suppose that $u$ is a nonnegative function which is harmonic with respect to
 $X$ in $D$. For any ball $B:=B(x_0, r)$ with $B\subset \overline{B}\subset D$,
the functions $u_n$, $n\ge 1$, defined by
$$
u_n(x):=\E_x\left[(u\wedge n)(X_{\tau_{
B}})\right], \qquad x\in \R^d\, ,
$$
are bounded functions which are harmonic with respect to $X$ in $B$.
Applying Theorem \ref{hi} to $v_n(x):=u(x)-u_n(x)=
\E_x\left[(u-(u\wedge n))(X_{\tau_{
B}})\right]$, it is easy to see that $u_n$
converges to $u$ uniformly in $B(x_0, r/2)$.
Thus $u$ is continuous in $D$.
This implies that all nonnegative functions that are harmonic in
$D$ with respect to $X$ are continuous.

A subset $D$ of $\R^d$ is said to be Greenian (for $X$) if $X^{D}$ is transient.
By \cite[Theorem 3.1]{CKK2} $X^D$ has  H\"older continuous transition densities $p_D(t,x,y)$.
For any Greenian open  set $D$ in $\R^d$, let $G_D(x,y)=\int_0^{\infty}p_D(t,x,y)dt$
be the Green function of $X^D$.
Then $G_D(x,y)$ is finite off the diagonal $D\times D$.
Furthermore, $x\mapsto G_D(x,y)$ is harmonic in $D\setminus \{y\}$
with respect to $X$ and therefore continuous.
Using the L\'{e}vy system for $X$, we know that for every
Greenian open subset $D$ and  every $f \ge 0$ and $x \in D$,
\begin{equation}\label{newls}
\E_x\left[f(X_{\tau_D});\,X_{\tau_D-} \not= X_{\tau_D}  \right]
=  \int_{\overline{D}^c} \int_{D}
G_D(x,z) J_X(z-y) dz f(y)dy.
\end{equation}
We define the Poisson kernel
\begin{equation}\label{PK}
K_D(x,y)\,:=  \int_{D}
G_D(x,z) J_X(z-y) dz, \qquad (x,y) \in
D \times  {\overline{D}}^c.
\end{equation}
Thus \eqref{newls} can be simply written as
$$
\E_x\left[f(X_{\tau_D});\,X_{\tau_D-} \not= X_{\tau_D}  \right]
=\int_{\overline{D}^c} K_D(x,y)f(y)dy.
$$

The following result will be used later in this paper.

\begin{lemma}
\label{l:KSV7prop4.7}
There exist $c_1=c_1(\Psi, \gamma_1, \gamma_2)>0$ and $c_2=c_2(\Psi, \gamma_1, \gamma_2)>0$
such that for every $r\in (0, 1]$
and $x_0\in \R^d$,
\begin{eqnarray}\label{e:KSV7-4.7}
K_{B(x_0, r)}(x, y)&\le& c_1j(|y-x_0|-r)(\phi(r^{-2})\phi((r-|x-x_0|)^{-2}))^{-1/2}\nonumber\\
&\le&c_1j(|y-x_0|-r)\phi(r^{-2})^{-1}
\end{eqnarray}
for all $(x, y)\in B(x_0, r)\times \overline{B(x_0, r)}^c$ and
\begin{equation}\label{e:KSV7-4.8}
K_{B(x_0, r)}(x_0, y)\ge c_2j(|y-x_0|)\phi(r^{-2})^{-1}, \qquad \mbox{for all } y
\in \overline{B(x_0, r)}^c.
\end{equation}
\end{lemma}

\pf This proof is exactly the same as that of \cite[Proposition 13.4.10]{KSV3}.
We provide the proof to show that only the monotonicity of $j$, \eqref{e:j-decay}, \eqref{e:doubling-condition}  and
Lemmas \ref{L3.2} and \ref{l:tau} are used.

Without loss of generality, we assume $x_0=0$. For $z\in B(0, r)$ and $r<|y|<2$,
$$
|y|-r\le |y|-|z|\le |y-z|\le |z|+|y|+r+|y|\le 2|y|,
$$
and for $z\in B(0, r)$ and $y\in B(0, 2)^c$,
$$
|y|-r\le |y|-|z|\le |y-z|\le |z|+|y|+r+|y|\le |y|+1.
$$
Thus by the monotonicity of $j$, \eqref{e:j-decay} and \eqref{e:doubling-condition}, there exists
a constant $c>0$ such that
$$
cj(|y|)\le j(|z-y|)\le j(|y|-r), \qquad (z, y)\in B(0, r)\times \overline{B(0, r)}^c.
$$
Applying the above inequalities, Lemmas \ref{L3.2} and \ref{l:tau} to \eqref{PK}, we immediately
get the assertion of the lemma.
\qed

As in \cite{KSV2}, to deal with $\kappa$-fat open set, we need the following form of Harnack inequality.

\begin{thm}\label{HP2}
Let $L>0$. There exists a positive constant $c=c(L, \Psi, \gamma_1, \gamma_2)>1$ such that
the following is true: If $x_{1}, x_{2}\in \R^d$ and
$r\in (0,1)$ are such
that $|x_{1}-x_{2}|< Lr$, then for every nonnegative function $u$ which
is harmonic with respect to $X$ in $B(x_{1}, r)\cup B(x_{2},r)$, we
have
$$
c^{-1}u(x_{2})\,\leq \,u(x_{1}) \,\leq\, c u(x_{2}).
$$
\end{thm}

\pf
Let $r \in (0,1]$,
$x_{1}, x_{2}\in \R^d$ be such that $|x_1-x_2|< Lr$
and let $u$ be a nonnegative function which is harmonic in
$B(x_{1}, r)\cup B(x_{2},r)$ with respect to $X$.
If $ |x_{1}- x_{2}|  <  \frac14 r$, then since $r<1$, the
theorem is true by Theorem \ref{hi}. Thus we only need to
consider the case when $\frac14 r \le |x_{1}- x_{2}| \le L r$ with
$L > \frac14$.

Let $w\in B(x_1, \frac{r}{8})$. Because $|x_2-w| \le
|x_1-x_2|+|w-x_1| < (L+\frac18)r \le 2Lr$, first using the monotonicity of
$j$ and \eqref{e:KSV7-4.8},
then using \eqref{e:asmpbofjat0} and Lemma \ref{l:phi-property}, we get
\begin{equation}\label{e:asdads1}
K_{B(x_2, \frac{r}8)}(x_2,w)  \,\ge\, c_1\,
j(2Lr)\phi(r^{-2})^{-1}\,  \,\ge\, c_2\,
r^{-d} \frac{\phi((2Lr)^{-2})}{\phi(r^{-2})} \ge c_3r^{-d}.
\end{equation}
For any $y\in B(x_1, \frac{r}8)$, $u$ is regular harmonic in $B(y,
\frac{7r}8)\cup B(x_1, \frac{7r}8)$. Since $|y-x_1|< \frac{r}8$, by
the already proven part of this theorem,
\begin{equation}\label{e:asdads2} u(y)\ge c_2
u(x_1), \quad y\in B(x_1, \frac{r}8),
\end{equation}
for some constant $c_2>0$. Therefore, by \eqref{newls} and
\eqref{e:asdads1}--\eqref{e:asdads2},
\begin{eqnarray*}
u(x_2) &=& \E_{x_2}\left[u(X_{\tau_{B(x_2,\frac{r}8) }})\right]
  \ge \E_{x_2}\left[u(X_{\tau_{B(x_2,\frac{r}8)}});
X_{\tau_{B(x_2,\frac{r}8)}} \in B(x_1,\frac{r}8) \right]\\
  &\ge &c_2 \,u(x_1)\, \P_{x_2}\left(X_{\tau_{B(x_2,\frac{r}8)}}
\in B(x_1,\frac{r}8) \right)
 = c_2\, u(x_1) \int_{B(x_1,\frac{r}8)} K_{B(x_2,\frac{r}8)} (x_2,w)\, dw \\
 &\ge  &c_3\, u(x_1) \left|B(x_1,\frac{r}8)\right|
r^{-d} \,= \, c_4\, u(x_1).
\end{eqnarray*}
Thus we have proved
the right-hand side inequality in the conclusion of the theorem. The
inequality on the left-hand side follows by symmetry.
\qed

For notational convenience,  we define
\begin{equation}\label{e:1.8}
\Phi(r)=\frac1{\phi(r^{-2})}, \quad r>0.
\end{equation}
The inverse function of $\Phi$ will be denoted by the usual notation $\Phi^{-1}(r)$.

Our process $X$ belongs to the class of Markov processes considered in  \cite{CKK2}. Thus
we have the following two-sided
estimates for $p(t, x, y)$ from
\cite{CKK2}. The proof is the same as that of \cite[Proposition 2.2]{CKS}.

\begin{prop}
\label{p:new1}
For any $T>0$, there exists $c_1=c_1(T, \Psi, \gamma_1, \gamma_2)>0$
such that
\begin{equation}\label{e:offdiag}
p(t, x, y) \leq
c_1 \, (\Phi^{-1}(t))^{-d}
\qquad  \hbox{ for all } (t, x, y)\in (0, T]\times \R^d \times \R^d.
\end{equation}
For any $T, R>0$, there exists
$c_2=c_2(T, R, \Psi, \gamma_1, \gamma_2)>1$ such that
for all $(t, x, y)\in [0, T]\times \R^d\times \R^d$ with $|x-y| <R,$
\begin{align}\label{stssbound}
c^{-1}_2\left( (\Phi^{-1}(t))^{-d}\wedge  \frac{t}{|x-y|^d \, \Phi (|x-y|)}
\right)
\le p(t, x, y) \le
c_2 \left((\Phi^{-1}(t))^{-d}\wedge \frac{t}{|x-y|^d \, \Phi (|x-y|)}\right).
\end{align}
\end{prop}

Our argument to obtain upper bound on the Green functions of bounded open sets is similar to those in \cite{CKS2, KS}. We give the
details here for the completeness.

\begin{lemma}\label{dec}
For every bounded open set $D$,
the Green function $G_D(x, y)$
is finite and continuous off the diagonal of $D\times D$ and
there exists $c =c(\mathrm{diam}(D), \Psi, \gamma_1, \gamma_2)\ge 1$
such that for all $x, y\in D$,
\begin{equation}\label{G_bd}
G_D(x, y) \, \le \, c \frac{\Phi(|x-y|)}{|x-y|^d}=\frac{c}{|x-y|^d\phi(|x-y|^{-2})}.
\end{equation}
\end{lemma}

\pf
Put $L:={\rm diam} (D)$. By \eqref{stssbound}, for
 every $x \in D$ we have
\begin{align*}
&\P_x ( \tau_D \le 1)
 \ge  \P_x ( X_1 \in  \bR^d \setminus D)
 =  \int_{ \bR^d \setminus D} p(1,x,y) dy\\
&\ge  c_1 \int_{ \bR^d \setminus D}
\left(1 \wedge \frac{1}{|x-y|^d \, \Phi (|x-y|)}\right) dy
 \ge  c_1 \int_{\{\,|z| \,\ge\, L \,\}}
\left(1 \wedge  \frac{1}{|z|^d \, \Phi (|z|)} \right)  dz=c_2>0.
\end{align*}
Thus
$$
\sup_{x \in D}\int_Dp_D(1, x, y)dy\,=\,
 \sup_{x \in D} \P_x ( \tau_D > 1) \,<\, 1.
$$
Now the Markov property of $X$  implies that there
exist positive constants $c_3$ and $c_4$
such that
$$
\int_Dp_D(t, x, y)dy\,\le\, c_3e^{-c_4t} \quad
\hbox{for all } (t, x)\in (0, \infty)\times D.
$$
Thus combining this, \eqref{stssbound} and the semigroup property, we have  that
for any $(t, x, y)\in (1, \infty)\times D\times D$,
\begin{eqnarray*}
p_D(t, x, y)&=&\int_Dp_D(t-1, x, z)p_D(1, z, y)dz
\le c_5\int_Dp_D(t-1, x, z)dz
\le c_6e^{-c_4t}.
\end{eqnarray*}
 This with \eqref{stssbound} implies that
 such that for any $(x, y)\in  D\times D$,
\begin{equation}\label{est:6.1}
G_D(x, y)=\int^{\infty}_0 p_D(t, x, y)dt \,\le\, c_7
\int_0^1\left((\Phi^{-1}(t))^{-d}\wedge  \frac{t}{|x-y|^d \, \Phi (|x-y|)}\right) dt +c_7.
\end{equation}
By the proof of \cite[Theorem 6.1]{CKK2},
$$
\int_0^1\left((\Phi^{-1}(t))^{-d}\wedge  \frac{t}{|x-y|^d \, \Phi (|x-y|)}\right) dt \asymp
 \frac{\Phi(|x-y|)}{|x-y|^d}.
$$
Therefore the Green function
$
G_D(x, y)=\int^{\infty}_0 p_D(t, x, y)dt
$
is finite and continuous off the diagonal of $D\times D$.
Furthermore, by \eqref{e:H1}, we have $\inf_{a \le {\rm diam}(D)} \Phi(a)a^{-d} >0$.
Consequently, \eqref{G_bd} holds.
\qed

For interior lower bound on the Green function, we use some recent results from \cite{CKS}.
The next result is an analog of \cite[Proposition 3.6]{CKS}, which is the main result of  \cite[Section 3]{CKS}.
Even though it is assumed in \cite{CKS} that $X$ is rotationally symmetric and its L\'evy density
satisfies a little stronger assumption than in this paper,
all the arguments of  \cite[Section 3]{CKS} only use
the results in \cite{CKK2},  \eqref{e:j-decay}, \eqref{e:psi1}, \eqref{e:asmpbofjat0}, \eqref{e:H1},
and the semigroup property.
Thus by following the same arguments line by line, one can prove the next proposition.
We omit the details.
We note in passing that $D$ is not necessarily bounded in the next proposition.

\begin{prop}\label{step31}
Let $T$ and $a$ be positive constants.
There exists $c=c(T, a, \Psi, \gamma_1, \gamma_2)>0$ such that for any open set $D$,
$$
p_D(t, x, y)\,\ge  \,c \, ((\Phi^{-1}(t))^{-d} \wedge  { t}{j(|x-y|)})
$$
for every $(t, x, y)\in
(0, T]\times D\times D$ with
$ \delta_D(x) \wedge \delta_D (y) \ge a \Phi^{-1}(t)$.
\end{prop}

\begin{lemma}\label{G:g3}
For every $L , T>0$, there exists $c=c(T, L, \Psi, \gamma_1, \gamma_2)>0$ such
that for any bounded open set $D$ with $\mathrm{diam}(D) \le T$,
$x, y\in D$ with $|x-y| \le L ( \delta_D(x) \wedge \delta_D(y))$,
\begin{equation}\label{e:g3}
G_D(x,y) \ge c \frac{\Phi(|x-y|)}{|x-y|^d}=\frac{c}{|x-y|^d\phi(|x-y|^{-2})}.
\end{equation}
\end{lemma}

\pf
Without loss of generality, we assume $L\ge 1$,
$\delta_D(y) \le \delta_D(x)$ and $\mathrm{diam}(D) \le T$.
By Proposition \ref{step31} and \eqref{e:asmpbofjat0},
there exists $c_1=c_1(T, \Psi, \gamma_1, \gamma_2)>0$
such that for all $(t, z, w)\in
(0, T]\times D\times D$ with
$ \delta_D(z) \wedge \delta_D (w) \ge  \Phi^{-1}(t)$,
\begin{align}
p_D(t, z, w)\,\ge  \,c_1 \, \left((\Phi^{-1}(t))^{-d} \wedge  \frac{ t}{ |z-w|^d \Phi( |z-w| )}\right). \label{e:estpd}
\end{align}
Using this we have
\begin{align}
&G_D(x,y) \ge \int_0^T p_D(t,x,y)dt \ge \int_0^{  \Phi(\delta_D(y))} p_D(t,x,y)dt\nn\\
& \ge c_1
\int_0^{ \Phi(\delta_D(y))} (\Phi^{-1}(t))^{-d} \wedge\frac{t}{|x-y|^d\Phi(|x-y|)}dt\nn\\
&\ge c_1
\int_0^{\Phi(L^{-1}|x-y|)} (\Phi^{-1}(t))^{-d} \wedge \frac{t}{|x-y|^d\Phi(|x-y|)}dt. \label{e:Glbd}
\end{align}
Let $r=|x-y|$. By the change of variable
$u= \frac{\Phi(r)}{t}$ and the fact that $t\to \Phi(t)$ is increasing,
\begin{align}
\int_0^{\Phi(L^{-1}r)} (\Phi^{-1}(t))^{-d} \wedge \frac{t}{r^d\Phi(r)}dt&=
\frac{\Phi(r)} {r^{d}}  \int_{\Phi(r)/\Phi(L^{-1}r)}^\infty
 u^{-2}\left(   \left( \frac{ r}{   \Phi^{-1} (  u^{-1} \Phi(r))}  \right)^d  \wedge u^{-1} \right)du \nn\\
 &=
 \frac{\Phi(r)} {r^{d}}  \int_{\Phi(r)/\Phi(L^{-1}r)}^\infty
 u^{-3}du.\label{e:Glbd1}
\end{align}
Since by \eqref{e:H1} and \eqref{e:1.8}
\begin{align}
 \int_{\Phi(r)/\Phi(L^{-1}r)}^\infty
 u^{-3}du \ge  \int_{c_2 L^{2 \delta_1}}^\infty u^{-3}du>0,
\end{align}
we conclude from \eqref{e:Glbd} and \eqref{e:Glbd1} that
\eqref{e:g3} holds.
\qed

It follows from \eqref{e:H1}, Lemmas \ref{G:g3} and \ref{l:tau} that, for all $r\in (0, 1)$ and all $\alpha\in (0, \pi)$,
there exists $c=c(\alpha)>0$ such that for all cones $V$ of angle $\alpha$ with
vertex at the origin,
\begin{align*}
&\E_0\int^{\tau_{B(0, r)}}_01_{V}(X_s)ds \ge  \int_{V \cap B(0, r/2)} G_{B(0, r)}(0, y)dy \ge
  c_1    \int_{V\cap B(0, r/2)}     \frac{1}{|y|^d\phi(|y|^{-2})}  dy \\
&
  \ge
  c_2    \int_{0}^{r/2}     \frac{1}{r\phi(r^{-2})}  dr
    \ge c_3 \phi(r^{-2})\ge c_4\E_0\tau_{B(0, r)}.
\end{align*}
Thus $X$ satisfies hypothesis {\bf H} in \cite{Sz00}. It follows from \eqref{e:psi1},
\eqref{e:asmpbofjat0} and Lemma \ref{L3.2} that for every cone $V$ with vertex at the origin,
$$
\int_{V\cap B(0, 1)}\E_0\tau_{B(0, |y|)}J_X(y)dy \ge c
\int_{V\cap B(0, 1)}
|y|^{-d}dy =
\infty.
$$
Therefore it follows from \cite[Theorem 1]{Sz00} that if $V$ is a Lipschitz open set in $\R^d$
and $D$ is an open subset of $V$,
$$
\P_x(X_{\tau_D}\in \partial V)=0, \qquad x\in D.
$$
Using Theorem \ref{hi},
 the proof of the next result is the same as that of
\cite[Proposition 13.4.11]{KSV3}. So we omit the proof.

\begin{prop}\label{p:Poisson2}
For every $a \in (0,1)$, there exists
$c=c(\Psi, \gamma_1, \gamma_2, a)>0$ such that for every
$r \in (0,1]$,
$x_0 \in {\mathbb R}^d$ and $x_1, x_2 \in B(x_0, ar)$,
$$
K_{B(x_0,r)}(x_1,y) \,\le\, c  K_{B(x_0,r)}(x_2,y), \qquad
\text{a.e. } y \in
\overline{B(x_0,r)}^{\, c}\, .
$$
\end{prop}

\begin{prop}\label{p:Poisson3}
For every $a \in (0,1)$, there exists
$c=c(\Psi, \gamma_1, \gamma_2, a)>0$ such that for every
$r \in (0, 1]$
and $x_0 \in \R^d$,
$$
K_{B(x_0,r)}(x,y)\le c\ r^{-d} \left(\frac{\phi((|y-x_0|-r)^{-2})}{\phi(r^{-2})}\right)^{1/2}
$$
for
all $x\in B(x_0,ar)$ and
a.e. $y$ such that $r<|x_0-y|<2r$.
\end{prop}

\pf
Using \eqref{e:asmpbofjat0}, Lemmas \ref{l:j-upper},  \ref{l:tau} and Proposition \ref{p:Poisson2},
the proof is exactly the same as that of \cite[Proposition 4.9]{KSV7}. We omit it.
\qed

\begin{lemma}\label{l2.1}
For every $a \in (0, 1)$, there exists a constant
$c=c(\Psi, \gamma_1, \gamma_2,  a)>0$ such that
for any $r\in (0, 1)$
and any open set $D$ with $D\subset B(0, r)$ we have
$$
{\P}_x\left(X_{\tau_D} \in B(0, r)^c\right) \,\le\, c\,
\phi(r^{-2})\int_D G_D(x,y)dy, \qquad x \in D\cap B(0, ar)\, .
$$
\end{lemma}
\pf
Using Lemma
\ref{l:lnew},  the proof is exactly the same as that of \cite[Lemma 4.10]{KSV7}. We omit it.
\qed

With these preparations in hand, we can repeat the argument of \cite[Section 5]{KSV7} to
get the following form of the boundary Harnack principle established in \cite{KSV7}.
We omit the details.
Note that the open set $D$ in the next result is not necessarily bounded.

\begin{thm}\label{ubhp}
There exists $C_1= C_1(\Psi, \gamma_1, \gamma_2)\ge 1$ such that
the following hold for all $r\in (0,1)$.
\begin{itemize}
    \item[(i)]
For every $z_0 \in \R^d$, every open set $U \subset B(z_0,r)$ and for any nonnegative function $u$ in $\RR^d$ which is regular harmonic in $U$ with respect to $X$ and vanishes a.e.~in $U^c \cap B(z_0, r)$ it holds that
\begin{eqnarray*}
C_1^{-1}\E_x[\tau_{U}] \int_{B(z_0, r/2)^c}  j(|y-z_0|)  u(y)dy
\,\le \,u(x)\,\le\,
C_1 \,\E_x[\tau_{U}] \int_{B(z_0, r/2)^c} j(|y-z_0|)  u(y)dy
\end{eqnarray*}
for every  $ x \in U \cap B(z_0, r/2)$.

    \item[(ii)] For every $z_0 \in \R^d$, every open set $D\subset \R^d$, every $r\in (0,1)$
    and for any nonnegative functions $u, v$ in $\R^d$ which are regular harmonic in $D\cap
    B(z_0, r)$ with respect to $X$ and vanish a.e. in $D^c \cap B(z_0, r)$, we have
        $$
        \frac{u(x)}{v(x)}\,\le C_1^4\,\frac{u(y)}{v(y)}, \qquad x, y\in D\cap B(z_0, r/2).
        $$

    \item[(iii)] For every $z_0 \in \R^d$, every
  Greenian open set
     $D\subset \R^d$, every $r\in (0,1)$, we have
    $$
    K_D(x_1, y_1)K_D(x_2, y_2) \le C_1^4 K_D(x_1, y_2)K_D(x_2, y_1)
    $$
    for all $x_1, x_2 \in D \cap B(z_0, r/2)$ and a.e. $y_1, y_2 \in
    \overline{D}^c \cap B(z_0, r)^c$.
\end{itemize}
  \end{thm}

In the next two results, we will assume that $X$ is transient
and will use $G(x, y)$ to denote the Green function of $X$.
Note that  $G(x,y)=G(y,x)$ by symmetry and $G(x,y)=G(0,y-x)$ by translation invariance.
Since we only assume that $X$ is symmetric, rather than unimodal, the next result does not follow from \cite{G}.

\begin{thm}\label{t:greenh}
For every $M \ge 1$  there exists a constant
$C_2(M)=C_2(M, \Psi, \gamma_1, \gamma_2)>0$
such that for all  $x\in B(0, M)$,
$$
C_2(M)^{-1} \frac{\Phi(|x|)}{|x|^{d}} \, \le\, G(x,0)
\,\le \, C_2(M)
\frac{\Phi(|x|)}{|x|^{d}}\, .
$$
\end{thm}
\pf
In this proof, we always assume that $x\in B(0, M)$.
It follows from Lemma \ref{G:g3}
that
$$
G(x,0) \ge G_{B(0, 2M)}(x,0) \ge  c_1 \frac{\Phi(|x|)}{|x|^{d}}.$$
On the other hand, by the strong Markov property and Lemma \ref{dec},
$$
G(x,0)=G_{B(0, 2M)}(x,0)+\E_x[G(X_{\tau_{B(0, 2M)}},0)] \le c_2
\frac{\Phi(|x|)}{|x|^{d}}
+\E_x[G(X_{\tau_{B(0, 2M)}},0)].
$$
Choose $x_0=(0, \dots, 0, M/2)$. By  Theorem \ref{HP2} and the strong Markov property, we have
$$
\E_x[G(X_{\tau_{B(0, 2M)}},0)]  \le c_3 \E_{x_0} [G(X_{\tau_{B(0, 2M)}},0)] \le
c_3 G(x_0,0)  \le c_4 < \infty.
$$
Furthermore, by \eqref{e:H1}, we have $\inf_{a \le 2M} \Phi(a)a^{-d} >0$. Thus
we conclude that
$$
G(x,0)=G_{B(0, 2M)}(x,0)+\E_x[G(X_{\tau_{B(0, 2M)}},0)] \le c_2
\frac{\Phi(|x|)}{|x|^{d}}
+c_4 \le c_5 \frac{\Phi(|x|)}{|x|^{d}}.
$$
\qed

The following result will be needed in the proof of Theorem \ref{t:m}.

\begin{lemma}\label{l:gfnatinfty}
If $X$ is transient, then
$$
\lim_{x\to\infty}G(x, 0)=0.
$$
\end{lemma}

\pf It follows from the transience assumption and
\eqref{e:psi3} that $1/\Psi_X$ is locally
integrable in $\R^d$, thus by \cite[Proposition 13.23]{BF}
that the semigroup of $X$ is integrable (in the sense of \cite{BF}).
Therefore it follows from \cite[Proposition 13.21]{BF} that, for any nonnegative
continuous function $f$ on $\R^d$ with compact support, $\lim_{x\to\infty}Gf(x)=0$.

By Theorem \ref{hi}, for any $x\in B(0, 4)^c$,
$$
G(x, 0)\,\le\, c\,G(x, y), \qquad y\in B(0, 2).
$$
Take a nonnegative function $f$ with support in $\overline{B(0, 2)}$ such that $f$ is identically 1 on $B(0, 1)$. Then
$$
G(x, 0)\int_{B(0, 2)}f(y)dy=\int_{B(0, 2)}G(x, 0)f(y)dy\le c
\int_{B(0, 2)}G(x, y)f(y)dy=Gf(x),
$$
which yields
$$
G(x, 0)\le c\frac{Gf(x)}{\int_{B(0, 2)}f(y)dy}.
$$
Therefore $\lim_{x\to\infty}G(x, 0)=0$.
\qed

\begin{remark}\label{r:h2}{\rm
Note that several results of this section are stated only for small radii $r$, namely $r\in (0,1]$.
This is, of course, a consequence of the scaling condition ({\bf H}) at infinity
which governs the behavior of the process for small time and small space.
If we want to study the large time and large space behavior of $X$,
we would need to add the following scaling condition on $\Psi$  near the origin too:

({\bf H2}):
There exist constants $0<\delta_3\le \delta_4 <1$ and $a_3, a_4>0$  such that
$$
a_3\lambda^{2\delta_4} \Psi(t) \le \Psi(\lambda t) \le a_4 \lambda^{2\delta_3} \phi(t),
\quad \lambda \le 1, t \le 1\, .
$$
One consequence of adding condition {\bf (H2)} to condition ({\bf H}) is that many results that were valid for small $r$ only will hold true for all $r>0$. For future reference, we list below precisely which statements are true.

First note that if both ({\bf H}) and ({\bf H2}) hold, there exist  $a_5, a_6>0$  such that
\begin{equation}\label{e:H2}
a_5 \left(\frac{R}{r}\right)^{2(\delta_1 \wedge \delta_3)} \le \frac{\Psi(R)}{\Psi(r)} \le a_6 \left(\frac{R}{r}\right)^{2(\delta_2 \vee \delta_4)}, \quad a>0,\ 0<r<R<\infty\, ,
\end{equation}
cf.~\cite[(2.6)]{KSV8}.

In the remainder of this remark, in addition to all conditions from Section \ref{sec:1}, namely \eqref{e:j-decay}, \eqref{e:psi1} and {\bf (H)}, we also assume ({\bf H2}). Then $X$ satisfies the assumptions in \cite{CK}. Furthermore, it follows from \cite[(15), Corollary 23, Proposition 28]{BGR1} and \eqref{e:H2} that  \eqref{e:asmpbofjat0} and \eqref{e:doubling-condition} hold for all $r>0$. That is, under both ({\bf H}) and ({\bf H2})
\begin{equation}\label{e:asmpbofjat0H2}
j(r)\asymp \frac{\phi(r^{-2})}{r^d}
\quad \hbox{for } r>0,
\end{equation}
and
\begin{equation}\label{e:doubling-conditionH2}
j(r)\le c j(2r)\, ,\quad r>0\, .
\end{equation}
The following results are now true for all $r>0$: Theorem \ref{hi}, Lemma \ref{l:KSV7prop4.7}, Theorem \ref{HP2}, Proposition  \ref{p:Poisson2}, Proposition \ref{p:Poisson3}, Lemma \ref{l2.1} and Theorem \ref{ubhp}. In order to prove Theorem \ref{hi}, we also have to use  \cite[Theorem 4.12]{CK}, while in the proof of Theorem \ref{HP2} we use \eqref{e:asmpbofjat0H2} instead of \eqref{e:asmpbofjat0}.
Furthermore, using results in \cite{CK}, Proposition \ref{step31} and Lemma \ref{G:g3} hold with $T=\infty$.
Proofs of all other listed results stay the same.
}
\end{remark}

\section{Martin Boundary}

In this section we will use the boundary Harnack principle (Theorem \ref{ubhp})
and follow the well-established route, see \cite{B, KSV1, KSV6, KSV9},
to study the Martin boundary of an open set $D$ with respect to $X$. The key ingredient
is the oscillation reduction technique used in the proof of Lemma \ref{l:oscillation-reduction}.

Throughout this section we assume that $D\subset \R^d$ is a Greenian open set,
$Q\in \overline{D}$ and $D$ is $\kappa$-fat
at $Q$ for some $\kappa\in (0, \frac12)$, that is, there is $R>0$ such that for all $r\in
(0,R]$, there is a ball  $B(A_r(Q) , \kappa r)  \subset D\cap B(Q,r)$.
Without loss of generality, we assume that $R \le 1/2$.
After Lemma \ref{l:oscillation-reduction}, we will further assume that $Q\in\partial D$.

\begin{lemma}\label{l:lemma3.1}
There exist
$C_3=C_3(\Psi, \kappa, \gamma_1, \gamma_2)>0$ and
$\xi =\xi(\Psi, \kappa, \gamma_1, \gamma_2) \in (0,1)$ such that for every
$r\in (0, R)$
and any non-negative function $h$ in $\R^d$ which is harmonic in $D\cap B(Q,r)$ it holds that
\begin{equation}\label{e:lemma3.1}
    h(A_r(Q))\le C_3 (\phi(r^{-2}))^{-1}
    \xi^{k}
    \phi((\kappa/2)^{-2k}r^{-2})
    h(A_{(\kappa/2)^{k} r}(Q))\, ,\quad k=0,1,2,\dots \, .
\end{equation}
\end{lemma}

\pf Without loss of generality, we may assume that $Q=0$.
Fix $r\in (0, R)$. For $k=0,1,2,\dots $,
let $\eta_k=(\kappa/2)^{k} r$,  $A_k=A_{\eta_k}(0)$ and $B_k=B(A_k,  \eta_{k+1})$. Note that the balls $B_k$ are pairwise disjoint. By harmonicity of $h$, for every $k=0,1,2 \dots $,
\begin{eqnarray*}
    h(A_k)
    = \E_{A_k}\left[h(X_{\tau_{B_k}})\right]
    \ge \sum_{l=0}^{k-1} \E_{A_k}\left[h(X_{\tau_{B_k}}):\, X_{\tau_{B_k}}\in B_l\right]
   = \sum_{l=0}^{k-1} \int_{B_l} K_{B_k}(A_k, z)h(z)\, dz\, .
\end{eqnarray*}
By Theorem \ref{HP2}, there exists
$c_1=c_1(\Psi, \gamma_1, \gamma_2)>0$
 such that for every $l=0,1,2,\dots$,
$
    h(z)\ge c_1 h(A_l)$ for all  $ z\in B_l$ .
Hence
$$
      \int_{B_l}K_{B_k}(A_k, z)h(z)\, dz \ge c_1 h(A_l) \int_{B_l}K_{B_k}(A_k, z)\, dz\, , \quad
      0\le l\le k-1\, .
$$
By \eqref{e:KSV7-4.8}, we have
$$
  \int_{B_l}  K_{B_k}(A_k,z)dz\ge c_2\phi(\eta_k^{-2})^{-1}\int_{B_l} j(|2(A_k-z)|)dz\, , \quad  0\le l\le k-1\, .
$$
For $l=0, 1, \cdots, k-1$ and $z\in B_l$,
it holds that $|z|\le \kappa(\kappa/2)^{l}r
+(\kappa/2)^{(l+1)}r
=(3\kappa/2)(\kappa/2)^lr$.
Since $ |A_k|\le \kappa  \eta_k  $, we have  that $|A_k-z|\le  |A_k|+|z|< 3\eta_l$.
Together with \eqref{e:asmpbofjat0} and Lemma \ref{l:phi-property}, this implies that
    $
    j(|2(A_k-z)|)\ge c_3 |\eta_l|^{-d}\phi(\eta_k^{-2}) $ for every  $z\in B_l$ and  $0\le l\le k-1$.
Therefore,
$$
    \int_{B_l}K_{B_k}(A_k, z)\, dz \ge c_4  |\eta_l|^{-d}\phi(\eta_k^{-2})\phi(\eta_k^{-2})^{-1} |{B_l}| \ge c_5
    \phi(\eta_l^{-2})/\phi(\eta_k^{-2})
    \, , \quad 0\le l\le k-1\, .
$$
Hence,
$$
    \phi(\eta_k^{-2}) h(A_k)\ge c_5 \sum_{l=0}^{k-1}\phi(\eta_l^{-2}) h(A_l)\, ,\quad \textrm{for all }k=1,2,\dots ,
$$
where $c_5=c_5(\Psi,\gamma_1, \gamma_2, \kappa)$.
Let $a_k:=\phi(\eta_k^{-2}) h(A_k)$ so that $a_k\ge c_5 \sum_{l=0}^{k-1} a_l$.
Using the identity
$1+  c_5 \sum_{l=0}^{k-2} (1+c_5)^{l}= (1+c_5)^{k-1}$ for $k \ge3$,
by induction it follows that $a_k\ge c_5 (1+c_5)^{k-1} a_0$.
Let $\xi:=(1+c_5)^{-1} \in (0,1)$ so that $$
 \phi(r^{-2}) h(A_0)=a_0\le
(1+c_5)c_5^{-1}\xi^{k} a_k=(1+c_5)c_5^{-1}\xi^{k} \phi(\eta_k^{-2}) h(A_k).
$$
Then,
$$
h(A_r(0))\le(1+c_5)c_5^{-1}(\phi(r^{-2}))^{-1}\xi^{k} \phi(\eta_k^{-2}) h( A_{(\kappa/2)^{k} r}(0))\, .
$$
\qed

\begin{lemma}\label{l:lemma3.2}
There exists $c=c(\Psi, \kappa, \gamma_1, \gamma_2)>0$ such that for every
$r\in (0, R)$ and
every non-negative function $u$ in $\R^d$ which is regular harmonic in
$D\cap B(Q, r)$ with respect to $X$,
$$
u(A_r(Q))\ge c
\phi( r^{-2})^{-1}
\int_{
B(Q, r)^c}j(|z-Q|)u(z)dz.
$$
\end{lemma}

\pf Without loss of generality, we may assume that $Q=0$.
Fix $r\in (0, R)$ and let $A:=A_r(0)$. Since $u$ is regular
harmonic in $D\cap B(0, (1-\kappa/2)r)$ with respect to $X$, we have
\begin{eqnarray*}
u(A)&\ge&\E_A\left[u(X_{\tau_{D\cap B(0, (1-\kappa/2)r)}}):
X_{\tau_{D\cap B(0, (1-\kappa/2)r)}}\in B(0, r)^c \right]\\
&=&\int_{B(0, r)^c}K_{D\cap B(0, (1-\kappa/2)r)}(A, z)u(z)dz\\
&=&\int_{B(0, r)^c}\int_{D\cap B(0, (1-\kappa/2)r)}G_{D\cap B(0, (1-\kappa/2)r)}(A, y)J_X(y-z)dyu(z)dz.
\end{eqnarray*}
Since $B(A, \kappa r/2)\subset D\cap B(0, (1-\kappa/2)r)$, by the domain monotonicity of Green functions,
$$
G_{D\cap B(0, (1-\kappa/2)r)}(A, y)\ge G_{B(A, \kappa r/2)}(A, y), \qquad y\in B(A, \kappa r/2).
$$
Thus
\begin{eqnarray*}
u(A)&\ge&\int_{B(0, r)^c}\int_{B(A, \kappa r/2)}G_{B(A, \kappa r/2)}(A, y)J_X(y-z)dyu(z)dz\\
&=&\int_{B(0, r)^c}K_{B(A, \kappa r/2)}(A, z)u(z)dz\\
&\ge&c_1\phi((\kappa r/2)^{-2})^{-1}\int_{B(0, r)^c}j(|A-z|)u(z)dz,
\end{eqnarray*}
where in the last line we used \eqref{e:KSV7-4.8}.
Note that $|A-z|\le 2|z|$ for $z\in A(0, r, 1)$ and $|A-z|\le |z|+1$ for $z\in B(0, 1)^c$. Hence
by \eqref{e:j-decay}, \eqref{e:doubling-condition} and Lemma
\ref{l:phi-property},
$$
u(A)\ge c_2\phi(r^{-2})^{-1}\int_{B(0, r)^c}j(|z|)u(z)dz.
$$
\qed

\begin{lemma}\label{l:lemma3.3}
There exist $C_4=C_4(\Psi, \kappa, \gamma_1, \gamma_2)\ge 1$ and
$\xi=\xi(\Psi, \kappa, \gamma_1, \gamma_2)\in (0, 1)$ such that for any $r\in
(0, 1)$, and any non-negative
function $u$ on $\R^d$ which is regular harmonic in $D\cap B(Q, r)$
with respect to $X$ and vanishes in $D^c\cap B(Q, r)$ we have
$$
\E_x\left[u(X_{\tau_{D\cap B_k}}): X_{\tau_{D\cap B_k}}\in B(Q, r)^c\right]
\le C_4\xi^k u(x), \qquad x\in D\cap B_k,
$$
where $B_k:=B(Q, (\kappa/2)^kr)$, $k=0, 1, 2, \cdots$.
\end{lemma}

\pf Without loss of generality, we may assume that $Q=0$. Fix
$r\in (0, R)$.
Let $\eta_k:=(\kappa/2)^kr$, $B_k:=B(0, \eta_k)$ and
$$
u_k(x):=\E_x\left[u(X_{\tau_{D\cap B_k}}): X_{\tau_{D\cap B_k}}\in B(0, r)^c\right],
\qquad x\in D\cap B_k.
$$
Note that for $x\in D\cap B_{k+1}$,
\begin{eqnarray*}
u_{k+1}(x)&=&\E_x\left[u(X_{\tau_{D\cap B_{k+1}}}): X_{\tau_{D\cap B_{k+1}}}\in B(0, r)^c\right]\\
&=&\E_x\left[u(X_{\tau_{D\cap B_{k+1}}}): \tau_{D\cap B_{k+1}}=\tau_{D\cap B_k},
X_{\tau_{D\cap B_{k+1}}}\in B(0, r)^c\right]\\
&=&\E_x\left[u(X_{\tau_{D\cap B_{k}}}): \tau_{D\cap B_{k+1}}=\tau_{D\cap B_k},
X_{\tau_{D\cap B_{k}}}\in B(0, r)^c\right]\\
&\le& \E_x\left[u(X_{\tau_{D\cap B_{k}}}):
X_{\tau_{D\cap B_{k}}}\in B(0, r)^c\right].
\end{eqnarray*}
Thus
\begin{equation}\label{e:5.4}
u_{k+1}(x)\le u_k(x), \qquad x\in D\cap B_{k+1}.
\end{equation}
Let $A_k:=A_{\eta_k}(0)$.
Similarly we have
\begin{eqnarray*}
u_k(A_k)&=&\E_{A_k}\left[u(X_{\tau_{D\cap B_k}}): X_{\tau_{D\cap B_k}}\in B(0, r)^c\right]\\
&\le& \E_{A_k}\left[u(X_{\tau_{B_k}}): X_{\tau_{D\cap B_k}}\in B(0, r)^c\right]\\
&\le&\int_{B(0, r)^c}K_{B_k}(A_k, z)u(z)dz.
\end{eqnarray*}
By \eqref{e:KSV7-4.7},  we have
$$
K_{B_k}(A_k, z)\le c_1j(|z|-\eta_k)\phi(\eta_k^{-2})^{-1},
\qquad z\in B(0, r)^c.
$$
Note that $|z|-\eta_k\ge |z|/2$ for $z\in A(0, r, 2)$ and $|z|-\eta_k\ge |z|-1$ for $z\ge 2$.
Thus by \eqref{e:j-decay}, \eqref{e:doubling-condition} and the monotonicity of $j$,
\begin{equation}\label{e:5.5}
u_k(A_k)\le c_2\phi(\eta_k^{-2})^{-1}\int_{B(0, r)^c}j(|z|)u(z)dz,
\qquad k=1, 2, \cdots.
\end{equation}
By Lemma \ref{l:lemma3.2}, we have
\begin{equation}\label{e:5.6}
u(A_0)\ge c_3\phi(
\eta_1^{-2})^{-1}\int_{B(0, r)^c}j(|z|)u(z)dz.
\end{equation}
Thus \eqref{e:5.5}--\eqref{e:5.6} imply that
$$
u_k(A_k)\le c_4\phi(
\eta_1^{-2})/\phi(\eta_k^{-2})u(A_0).
$$
On the other hand, using Lemma \ref{l:lemma3.1}, we get
$$
u(A_0)\le c_5(\phi(r^{-2}))^{-1}\xi^{k} \phi(\eta_k^{-2})u(A_k).
$$
Combining the last two displays and using Lemma \ref{l:phi-property}, we get
$$
u_k(A_k)\le c_6\xi^{k} \frac{\phi(
\eta_1^{-2})}{\phi(r^{-2})}u(A_k)
\le c_7 \xi^k
u(A_k).
$$
By the boundary Harnack principle, Theorem \ref{ubhp} (i), we have
$$
\frac{u_k(x)}{u(x)}\le \frac{u_{k-1}(x)}{u(x)}\le c_8
\frac{u_{k-1}(A_{k-1})}{u(A_{k-1})}\le
 c_9\xi^k
$$
for $k=1, 2, \cdots$.
The proof is now complete.
\qed

Using the boundary Harnack principle and Lemma \ref{l:lemma3.3} (instead
of Lemmas 13 and 14 in \cite{B}), we can repeat the argument
of \cite[Lemma 16]{B} (which dealt with isotropic stable process) to get the following result.
We include the proof from \cite{B} to show that it does not depend on scaling property of
stable processes, and on the way make some constants explicit and provide the detailed
computation in the end of the induction argument.

\begin{lemma}\label{l:oscillation-reduction}
There exist $c=c(\Psi, \gamma_1, \gamma_2, \kappa)>0$ and
$\beta=\beta(\Psi, \gamma_1, \gamma_2, \kappa)>0$
such that for all
$r \in (0, R/2)$ and non-negative functions $u$ and $v$ on $\R^d$ which
are regular harmonic in $D\cap B(Q, 2r)$ with respect to $X$, vanish on
$D^c\cap B(Q, 2r)$ and satisfy $u(A_r(Q))=v(A_r(Q))$, the limit
$g=\lim_{D\ni x\to Q}\frac{u(x)}{v(x)}$ exists, and we have
\begin{equation}\label{e:zoran5.6}
\left|\frac{u(x)}{v(x)}-g \right|\le c\left(\frac{|x-Q|}{r}\right)^{\beta}, \qquad x\in D\cap B(Q, r).
\end{equation}
\end{lemma}

\pf Fix $r\in (0, R/2)$. Without loss of generality, we assume that $u(A_r(Q))=v(A_r(Q))=1$.
Let $D_0=D\cap B(Q,r)$. We start the proof by fixing several constants. We first choose
$c_1=c_1(\Psi, \gamma_1, \gamma_2)\ge 10$ such that
\begin{equation}\label{e:zoran5.13z}
\sup_{x\in D_0}\frac{u(x)}{v(x)}\le \left(1+\frac{c_1}{2}\right)\inf_{x\in D_0}\frac{u(x)}{v(x)}\, .
\end{equation}
This is possible because of Theorem \ref{ubhp} (it suffices to choose $1+c_1/2\ge C_1^2$). Let
$\delta=\delta(\Psi,\gamma_1, \gamma_2):=1-\frac12C^{-1}_1\in (1/2, 1)$,
where $C_1$ is the constant from Theorem \ref{ubhp}.
We further define $\epsilon=\epsilon(\Psi, \gamma_1, \gamma_2)\in (0,1/4)$ by
$$
\epsilon:=\frac{1-\delta}{20c_1}\, ,
$$
and choose $k_0=k_0(\Psi, \gamma_1, \gamma_2, \kappa)\in \N$ large enough so that
$(1-\xi^{k_0})^{-1}\le 2$ and $C_4 \xi^{k_0} \le \eps (4C_1)^{-1}$
where $C_4\ge 1$ and $\xi\in (0,1)$ are the constants from Lemma \ref{l:lemma3.3}.

For $k=0, 1, \dots$, define
$$
r_k=(\kappa/2)^{k_0k} r,
 \ \ B_k=B(Q, r_k), \ \ D_k=D\cap B_k, \ \
\Pi_k=D_k\setminus D_{k+1}, \ \ \Pi_{-1}=B_0^c.
$$
For $l=-1, 0, 1, \dots, k-1$, let
\begin{eqnarray}
u^l_k(x)&:=&\E_x\left[u(X_{\tau_{D_k}}); X_{\tau_{D_k}}\in \Pi_l\right], \qquad x\in \R^d,\label{e:zoran5.7}\\
v^l_k(x)&:=&\E_x\left[v(X_{\tau_{D_k}}); X_{\tau_{D_k}}\in \Pi_l\right], \qquad x\in \R^d.\label{e:zoran5.8}
\end{eqnarray}
Apply Lemma \ref{l:lemma3.3} with $\tilde{r}=r_{l+1}$ instead of $r$. Then $r_k=(\kappa/2)^{k_0 k}
r=(\kappa/2)^{k_0(k-l-1)}\tilde{r}$, hence for $k=0,1,2,\dots$ and $x\in D_k$
\begin{equation}\label{e:zoran5.9}
u^l_k(x)\le C_4(\xi^{k_0})^{k-l-1}u(x), \qquad l=-1, 0, 1, \dots, k-2.
\end{equation}
Then since
$\sum^{k-2}_{l=-1}(\xi^{k_0})^{k-1-l}=\xi^{k_0}\sum^{k-1}_{n=0}(\xi^{k_0})^n
\le \xi^{k_0}(1-\xi^{k_0})^{-1}\le 2\xi^{k_0}$, we have that for $k=1, 2, \dots$
and $x\in D_k$,
\begin{equation}\label{e:zoran5.10}
\sum^{k-2}_{l=-1}u^l_k(x)\le 2C_4\xi^{k_0}u(x).
\end{equation}
Since $2C_4 \xi^{k_0} \le  \eps/2 < 1/2$, for $k=1, 2, \dots, l=-1, 0, 1, \dots, k-2,$
we have $\sum^{k-2}_{l=-1}u^l_k(x) \le u_k^{k-1}(x)$ for all $x\in D_k$ and, by \eqref{e:zoran5.10}
\begin{equation}\label{e:zoran5.11}
u^l_k(x)\le \epsilon^{k-1-l}u^{k-1}_k(x), \qquad x\in D_k.
\end{equation}
By symmetry we also have that for $k=1, 2, \dots, l=-1, 0, 1, \dots, k-2,$ \begin{equation}\label{e:zoran5.12}
v^l_k(x)\le \epsilon^{k-1-l}v^{k-1}_k(x), \qquad x\in D_k.
\end{equation}

Define $\zeta=\zeta(\Psi, \gamma_1, \gamma_2)\in (1/2, 1)$ by
$$
\zeta=\sqrt{\frac{1+\delta}{2}}\, .
$$
We claim that for all $l=0, 1, \dots$
\begin{equation}\label{e:zoran5.13}
\sup_{x\in D_l}\frac{u(x)}{v(x)}\le (1+c_1\zeta^l)\inf_{x\in D_l}\frac{u(x)}{v(x)}.
\end{equation}
By Theorem \ref{ubhp} with $y=A_r(Q)$ we have that
\begin{equation}\label{e:zoran5.14}
C_1^{-1}\le \frac{u(x)}{v(x)}\le C_1, \qquad x\in D_0=D\cap B(Q, r).
\end{equation}
We show that \eqref{e:zoran5.13} and \eqref{e:zoran5.14} imply the statement of the lemma.
Indeed,
$$
\sup_{x\in D_l}\frac{u(x)}{v(x)}-\inf_{x\in D_l}\frac{u(x)}{v(x)}\le c_1\zeta^l
\inf_{x\in D_l}\frac{u(x)}{v(x)}\le C_1c_1\zeta^l, \qquad l=0, 1, \dots.
$$
Since the right-hand side goes to zero as $l\to\infty$, the same is valid for
the left-hand side proving that the limit $g=\lim_{D\ni x\to Q}\frac{u(x)}{v(x)}$ exists.
Further, for $x\in D\cap B(Q, r)$ there is a unique $l\ge 1$ such that $x\in \Pi_{l-1}\subset D_{l-1}$. Thus $x\notin D_l$
implying $|x-Q|\ge r_l=(\kappa/2)^{k_0l}r$,
that is
$$
l\ge \frac{\log\frac{r}{|x-Q|}}{k_0\log\frac{2}{\kappa}}.
$$

Let $\beta:=-\frac{\log\zeta}{k_0\log\frac{2}{\kappa}}$.
Then $\beta=\beta(\Psi, \gamma_1, \gamma_2, \kappa)$ and

$$
\left|\frac{u(x)}{v(x)}-g\right|
\le C_1c_1\zeta^l
\le C_1c_1  \zeta^{\frac{\log\frac{r}{|x-Q|}}{k_0\log\frac{2}{\kappa}}}
=C_1c_1 \left(\frac{r}{|x-Q|}\right)^{\frac{\log\zeta}{k_0\log\frac{2}{\kappa}}}
= C_1c_1\left(\frac{|x-Q|}r\right)^{\beta}.
$$

We now prove \eqref{e:zoran5.13} by induction.
By \eqref{e:zoran5.13z} we wee that \eqref{e:zoran5.13} holds for $l=0$.
Again by \eqref{e:zoran5.13z} and the fact that $\zeta>1/2$ we have
$$
\sup_{x\in D_1}\frac{u(x)}{v(x)}\le \sup_{x\in D_0}\frac{u(x)}{v(x)} \le \left(1+\frac{c_1}{2}\right)\inf_{x\in D_0}\frac{u(x)}{v(x)} \le \left(1+c_1 \zeta\right)\inf_{x\in D_0}\frac{u(x)}{v(x)}\le \left(1+c_1 \zeta\right)\inf_{x\in D_1}\frac{u(x)}{v(x)}\, ,
$$
hence \eqref{e:zoran5.13} holds also for $l=1$.
Let $k=0,1, \dots$ and assume that \eqref{e:zoran5.13} holds for $l=0, 1, 2, \dots, k$.
By the definitions \eqref{e:zoran5.7}--\eqref{e:zoran5.8}, and the regular harmonicity of $u$ and $v$,
we have
\begin{eqnarray}
u(x)&=&\sum^k_{l=-1}u^l_{k+1}(x), \qquad x\in
D_{k+1},\label{e:zoran5.15}\\
v(x)&=&\sum^k_{l=-1}v^l_{k+1}(x), \qquad x\in
D_{k+1}.\label{e:zoran5.16}
\end{eqnarray}
For any function $f$ on a set $A$ we define
$$
{\rm Osc}_Af=\sup_{x\in A}f(x)-\inf_{x\in A}f(x).
$$
Let $g(x):=u^k_{k+1}(x)/v^k_{k+1}(x)$, $x\in D_k$. We claim that
\begin{equation}\label{e:zoran5.17}
{\rm Osc}_{D_{k+2}}g\le \delta\, {\rm Osc}_{D_k}g.
\end{equation}
Recall that $\delta=1-\frac12C^{-1}_1\in (1/2, 1)$.
Let $m_1:=\inf_{x\in D_k}g(x)$ and $m_2:=\sup_{x\in D_k}g(x)$.
By  \eqref{e:zoran5.7}, \eqref{e:zoran5.8} and Theorem \ref{ubhp}, it holds that $0<m_1
\le m_2<\infty$.
If $m_1=m_2$, then
both ${\rm Osc}_{D_{k+2}}g$ and ${\rm Osc}_{D_k}g$ are zero
so \eqref{e:zoran5.17} holds trivially.
Otherwise, let
$$
\widetilde{g}(x):=\frac{g(x)-m_1}{m_2-m_1}=\frac{u^k_{k+1}(x)-m_1v^k_{k+1}(x)}
{(m_2-m_1)v^k_{k+1}(x)}, \qquad x\in D_k.
$$
Note that $\widetilde{g}$ is the quotient of two non-negative functions
regular harmonic in $D_{k+1}$ with respect to $X$. Clearly ${\rm Osc}_{D_k}\widetilde{g}=1$.
Furthermore,
\begin{equation}\label{e:zoran5.18}
{\rm Osc}_{D_{k+2}}g={\rm Osc}_{D_{k+2}}\widetilde{g}\cdot{\rm Osc}_{D_{k}}g.
\end{equation}
This is clear from $g(x)=(m_2-m_1)\widetilde{g}(x)+m_1=({\rm Osc}_{D_{k}}g)\widetilde{g}(x) +m_1$.
If $\sup_{D_{k+2}}\widetilde{g}(x)\le 1/2$, then ${\rm Osc}_{D_{k+2}}\widetilde{g}\le 1/2$, and it follows
from \eqref{e:zoran5.18} that
\begin{equation}\label{e:zoran5.19}
{\rm Osc}_{D_{k+2}}g\le \frac12 {\rm Osc}_{D_{k}}g.
\end{equation}
If, on the other hand, $\sup_{D_{k+2}}\widetilde{g}(x)> 1/2$, we apply Theorem \ref{ubhp}
to the functions $\widetilde{u}(x)=u^k_{k+1}(x)-m_1v^k_{k+1}(x)$ and $\widetilde{v}(x)=v^k_{k+1}(x)$
to conclude that
$$
C^{-1}_1\frac{\widetilde{u}(y)}{\widetilde{v}(y)}\le \frac{\widetilde{u}(x)}{\widetilde{v}(x)}
\le C_1\frac{\widetilde{u}(y)}{\widetilde{v}(y)}, \qquad x, y\in D_{k+1}.
$$
This can be written as
$$
C^{-1}_1\widetilde{g}(y)\le \widetilde{g}(x)\le C_1\widetilde{g}(y), \qquad x, y\in D_{k+1}.
$$
Hence, for all $x\in D_{k+2}$, we have
$$
\widetilde{g}(x)\ge C^{-1}_1\sup_{y\in D_{k+2}}\widetilde{g}(y)\ge\frac12C^{-1}_1.
$$
Therefore, $\inf_{y\in D_{k+2}}\widetilde{g}(y)\ge\frac12C^{-1}_1$, and since $\widetilde{g}\le 1$,
we get that
\begin{equation}\label{e:zoran5.20}
{\rm Osc}_{D_{k+2}}\widetilde{g}\le 1-\frac12C^{-1}_1=\delta.
\end{equation}
By \eqref{e:zoran5.18}--\eqref{e:zoran5.20} we get \eqref{e:zoran5.17}.

We claim that
\begin{equation}\label{e:zoran5.21}
\inf_{x\in D_{k+2}}\frac{u^k_{k+1}(x)}{v^k_{k+1}(x)}\ge \inf_{x\in D_{k+1}}\frac{u^k_{k+1}(x)}{v^k_{k+1}(x)}
\ge \inf_{x\in D_k}\frac{u(x)}{v(x)}\ge \inf_{x\in D_i}\frac{u(x)}{v(x)}, \quad i=0, \dots k.
\end{equation}
Indeed, let $\eta:= \inf_{x\in D_k}\frac{u(x)}{v(x)}$ so that $u(x)\ge \eta v(x)$ on $D_k$.
Then for $x\in D_{k+1}$,
$$
u^k_{k+1}(x)=\E_x[u(X_{\tau_{D_{k+1}}}); X_{\tau_{D_{k+1}}}\in \Pi_k]\ge
\eta \E_x[v(X_{\tau_{D_{k+1}}}); X_{\tau_{D_{k+1}}}\in \Pi_k]=\eta v^k_{k+1}(x),
$$
so the second inequality of \eqref{e:zoran5.21} holds.
The first and third equalities of \eqref{e:zoran5.21} are trivial.
Similarly we have
\begin{equation}\label{e:zoran5.22}
\sup_{x\in D_{k+2}}\frac{u^k_{k+1}(x)}{v^k_{k+1}(x)}\le  \sup_{x\in D_{k+1}}
\frac{u^k_{k+1}(x)}{v^k_{k+1}(x)}\le \sup_{x\in D_k}\frac{u(x)}{v(x)}.
\end{equation}
Combining \eqref{e:zoran5.21}, \eqref{e:zoran5.22} and \eqref{e:zoran5.17} we get
\begin{equation}\label{e:zoran5.23}
\frac{\sup_{x\in D_{k+2}}\frac{u^k_{k+1}(x)}{v^k_{k+1}(x)}}{\inf_{x\in D_{k+2}}\frac{u^k_{k+1}(x)}{v^k_{k+1}(x)}}
-1\le \delta\left(\frac{\sup_{x\in D_k}\frac{u(x)}{v(x)}}{\inf_{x\in D_k}\frac{u(x)}{v(x)}}-1\right).
\end{equation}
Using \eqref{e:zoran5.23} and \eqref{e:zoran5.13} with $l=k$ (the induction hypothesis) we obtain
\begin{equation}\label{e:zoran5.24}
\sup_{x\in D_{k+2}}\frac{u^k_{k+1}(x)}{v^k_{k+1}(x)}=(1+c_1\delta\rho\zeta^k)
\inf_{x\in D_{k+2}}\frac{u^k_{k+1}(x)}{v^k_{k+1}(x)}
\end{equation}
with a suitably chosen $\rho\in [0, 1]$ (independent of $x\in D_{k+2}$). Indeed,
$$
M:=\frac{\sup_{x\in D_{k+2}}\frac{u^k_{k+1}(x)}{v^k_{k+1}(x)}}{\inf_{x\in D_{k+2}}\frac{u^k_{k+1}(x)}{v^k_{k+1}(x)}}
\le 1+\delta\left(\frac{\sup_{x\in D_k}\frac{u(x)}{v(x)}}{\inf_{x\in D_k}\frac{u(x)}{v(x)}}-1\right)
\le 1+c_1\delta\zeta^k.
$$
Thus, $1\le M\le1+c_1\delta\zeta^k$ which implies that there exists $\rho=\rho(
c_1, \zeta, k)\in [0, 1]$ such that
$M=1+c_1\delta\rho\zeta^k$. Next, by symmetry and \eqref{e:zoran5.24} we have
\begin{equation}\label{e:zoran5.25}
\sup_{x\in D_{k+2}}\frac{v^k_{k+1}(x)}{u^k_{k+1}(x)}=(1+c_1\delta\rho\zeta^k)
\inf_{x\in D_{k+2}}\frac{v^k_{k+1}(x)}{u^k_{k+1}(x)}.
\end{equation}
Combining \eqref{e:zoran5.21} and \eqref{e:zoran5.13} we get
\begin{equation}\label{e:zoran5.26}
\sup_{x\in D_l}\frac{u(x)}{v(x)}\le (1+c_1\zeta^l)\inf_{x\in D_{k+2}}\frac{u^k_{k+1}(x)}{v^k_{k+1}(x)},
\qquad l=0, 1, \dots, k.
\end{equation}

Now we fix $x\in D_{k+2}$. Then by \eqref{e:zoran5.15} and \eqref{e:zoran5.16},
\begin{eqnarray}
\frac{u(x)}{v(x)}&=&\frac{\sum^k_{l=0}u^l_{k+1}(x)+u^{-1}_{k+1}(x)}{\sum^k_{l=0}v^l_{k+1}(x)+v^{-1}_{k+1}(x)}
\nonumber\\
&\le&\frac{\sum^k_{l=0}u^l_{k+1}(x)+\epsilon^{k+1}u^{k}_{k+1}(x)}{\sum^k_{l=0}v^l_{k+1}(x)}\nonumber\\
&\le& (1+\epsilon^{k+1})\frac{\sum^k_{l=0}u^l_{k+1}(x)}{\sum^k_{l=0}v^l_{k+1}(x)}, \label{e:zoran5.27}
\end{eqnarray}
where in the first inequality  we used \eqref{e:zoran5.11} with $l=-1$. Now we apply \cite[Lemma 15]{B}
with
\begin{eqnarray*}
&&u_0=u^k_{k+1}(x), \ \ v_0=v^k_{k+1}(x), \ \ a=(1+c_1\delta\rho\zeta^k)\inf_{y\in D_{k+2}}\frac{u^k_{k+1}(y)}
{v^k_{k+1}(y)},\\
&&u_i=u^{k-i}_{k+1}(x), \ \ v_i=v^{k-i}_{k+1}(x),\ \ b_i=(1+c_1\zeta^{k-i})
\inf_{y\in D_{k+2}}\frac{u^k_{k+1}(y)}{v^k_{k+1}(y)},\ \ \epsilon_i=\epsilon^i, \ \ i=1, \dots, k.\\
\end{eqnarray*}
We need to check the conditions of \cite[Lemma 15]{B}. \cite[(5.16)]{B}, $a\le b_i$, $i=1, \dots, k$, is
immediate. \cite[(5.18)]{B}, $u_0\le av_0$, is true by \eqref{e:zoran5.24}. \cite[(5.17)]{B} consists
of two parts: $u_i\le b_iv_i$ and $v_i\le \epsilon_iv_0$ for $i=1, \dots, k$.  The first part follows from
$$
\sup_{y\in D_{k+1}}\frac{u^{k-i}_{k+1}(y)}{v^{k-i}_{k+1}(y)}\le \sup_{y\in D_{k-i}}\frac{u(y)}{v(y)}
\le (1+c_1\zeta^{k-i})\inf_{y\in D_{k+2}}\frac{u^{k-i}_{k+1}(y)}{v^{k-i}_{k+1}(y)},
$$
where the first inequality is
\eqref{e:zoran5.22}
and the second inequality is
\eqref{e:zoran5.26}. The second part of \cite[(5.17)]{B} is precisely \eqref{e:zoran5.12}. Now
we apply \cite[Lemma 15]{B} to conclude that
\begin{eqnarray*}
\lefteqn{\sum^k_{i=0}u^{k-i}_{k+1}(x)}\\
&\le & \left[(1+c_1\delta\rho\zeta^k)\inf_{y\in D_{k+2}}\frac{u^{k}_{k+1}(y)}{v^{k}_{k+1}(y)}
+\sum^k_{i=1}c_1(\zeta^{k-i}-\delta\rho\zeta^k)\inf_{y\in D_{k+2}}\frac{u^{k}_{k+1}(y)}{v^{k}_{k+1}(y)}
\cdot\epsilon^i\right]\cdot\sum^k_{i=0}v^{k-i}_{k+1}(x)\\
&\le &\left[1+c_1\delta\rho\zeta^k+ c_1\sum^k_{i=1}\zeta^{k-i}\epsilon^i\right]\inf_{y\in D_{k+2}}\frac{u^{k}_{k+1}(y)}{v^{k}_{k+1}(y)}\cdot\sum^k_{i=0}v^{k-i}_{k+1}(x)\, .
\end{eqnarray*}
Hence, letting $\tau:= (1+\epsilon^{k+1})(1+c_1\delta\rho\zeta^k+
c_1\sum^k_{i=1}\zeta^{k-i}\epsilon^i)$ and applying \eqref{e:zoran5.27} we get
\begin{equation}\label{e:zoran5.28}
\frac{u(x)}{v(x)}\le (1+\epsilon^{k+1})\frac{\sum^k_{i=0}u^{k-i}_{k+1}(x)}{\sum^k_{i=0}v^{k-i}_{k+1}(x)}
\le \tau\cdot \inf_{y\in D_{k+2}}\frac{u^k_{k+1}(y)}{v^k_{k+1}(y)}.
\end{equation}
By symmetry,
\begin{equation}\label{e:zoran5.29}
\frac{v(x)}{u(x)}\le\tau\cdot \inf_{y\in D_{k+2}}\frac{v^k_{k+1}(y)}{u^k_{k+1}(y)}.
\end{equation}
Now, \eqref{e:zoran5.28}, \eqref{e:zoran5.24} and \eqref{e:zoran5.29}
imply that
\begin{eqnarray*}
\sup_{y\in D_{k+2}}\frac{u(y)}{v(y)}&\le &\tau\cdot \inf_{y\in D_{k+2}}\frac{u^k_{k+1}(y)}{v^k_{k+1}(y)}\\
&=&\frac{\tau}{1+c_1\delta\rho\zeta^k}\sup_{y\in D_{k+2}}\frac{u^k_{k+1}(y)}{v^k_{k+1}(y)}\\
&\le&\frac{\tau^2}{1+ c_1\delta\rho\zeta^k}\inf_{y\in D_{k+2}}\frac{u(y)}{v(y)}.
\end{eqnarray*}

The proof of the induction step will be finished
by showing that
\begin{equation}\label{e:ksv1}
\frac{\tau^2}{1+ c_1\delta\rho\zeta^k}\le 1+c_1\zeta^{k+2}.
\end{equation}
First note that
$$
\tau=(1+\epsilon^{k+1})\left(1+c_1\left(\delta\rho +\frac{\epsilon}{\zeta-\epsilon}\left(1-\left(\frac{\epsilon}{\zeta}\right)^k\right)\right)\zeta^k\right).
$$
Let
$$
\wt{\tau}:=(1+\epsilon^{k+1})\left(1+c_1\left(\delta\rho +\frac{\epsilon}{\zeta-\epsilon}\right)\zeta^k\right)\, ,
$$
so that $\tau\le \wt{\tau}$. Then
\begin{eqnarray*}
\tau^2\le\wt{\tau}^2 &=&(1+\epsilon^{k+1})^2\left(1+2c_1\left(\delta\rho+\frac{\epsilon}{\zeta-\epsilon}\right)\zeta^k +c_1^2\left(\delta\rho+\frac{\epsilon}{\zeta-\epsilon}\right)^2\zeta^{2k}\right)\\
&\le &(1+3\epsilon^{k+1})A\, ,
\end{eqnarray*}
where
$$
A:=1+c_1\left(2\delta\rho+\frac{2\epsilon}{\zeta-\epsilon}\right)\zeta^k +c_1^2\left(\delta\rho+\frac{\epsilon}{\zeta-\epsilon}\right)^2\zeta^{2k}\, .
$$
Further, let
$$
B:=(1+c_1\delta\rho\zeta^k)(1+c_1\zeta^{k+2})=1+c_1(\delta\rho+\zeta^2)\zeta^k+c_1^2(\delta\rho\zeta^2)\zeta^{2k}\, .
$$
We will show that $B-(1+3\epsilon^{k+1})A=B-A-3\epsilon^{k+1}A\ge 0$.
Since $ B-(1+3\epsilon^{k+1})A\le B-\tau^2$,
this will prove what want.
We now prove $B-(1+3\epsilon^{k+1})A=B-A-3\epsilon^{k+1}A\ge 0$.
Note that $\epsilon\in
(0, 1/200)\subset (0,1/4)$ and thus $\zeta-\epsilon>1/2-1/4=1/4$.
We will need the following estimate:
\begin{equation}\label{e:epsilon-zeta}
\frac{\epsilon}{\zeta-\epsilon}\le 4\epsilon <1\, .
\end{equation}
Since $\delta\in (1/2,1)$, $\rho\in [0,1]$ and $\zeta<1$ we have
\begin{equation}\label{e:upper-A}
A\le 1+c_1(2+2)+c_1^2(1+1)^2=1+4c_1+4c_1^2=(1+2c_1)^2\le (3c_1)^2=9c_1^2\, .
\end{equation}
Note that it follows from \eqref{e:epsilon-zeta} that
$$
\zeta^2-\delta-8\epsilon=\frac{1+\delta}{2}-\delta-\frac{2(1-\delta)}{5c_1}
=(1-\delta)\left(\frac12-\frac{2}{5c_1}\right)\ge \frac25(1-\delta)\, ,
$$
where we used that $1/2-2/(5c_1)\ge 1/2-2/50=23/50\ge 2/5$ which follows immediately from $c_1\ge 10$.
Thus,
\begin{eqnarray}\label{e:lower-B-A}
B-A&=&c_1\left(\delta\rho+\zeta^2 -2\delta\rho-\frac{2\epsilon}{\zeta-\epsilon}\right)\zeta^k +c_1^2\left(\delta\rho\zeta^2-\delta^2\rho^2-2\delta\rho\frac{\epsilon}{\zeta-\epsilon} -\frac{\epsilon^2}{(\zeta-\epsilon)^2}\right)\zeta^{2k}\nonumber \\
&=&c_1\left(\zeta^2 -\delta\rho-\frac{2\epsilon}{\zeta-\epsilon}\right)\zeta^k +c_1^2\left[\delta\rho\left(\zeta^2 -\delta\rho-\frac{2\epsilon}{\zeta-\epsilon}\right)-\frac{\epsilon^2}{(\zeta-\epsilon)^2}\right]\zeta^{2k}\nonumber\\
&>&c_1(\zeta^2-\delta\rho-8\epsilon)\zeta^k-c_1^2 \frac{\epsilon^2}{(\zeta-\epsilon)^2}\zeta^{2k}\nonumber\\
&\ge & c_1(\zeta^2-\delta-8\epsilon)\zeta^k-c_1^2(4\epsilon)^2 \zeta^{2k}\nonumber\\
&\ge &\frac{2}{5}c_1 (1-\delta)\left(\frac{1+\delta}{2}\right)^{k/2} -c_1^2\left(\frac{1-\delta}{5c_1}\right)^2\left(\frac{1+\delta}{2}\right)^{k}\nonumber\\
&=&\frac{2}{5}c_1 (1-\delta)\left(\frac{1+\delta}{2}\right)^{k/2} -\frac{(1-\delta)^2}{25}\left(\frac{1+\delta}{2}\right)^{k}\nonumber\\
&=&(1-\delta)\left(\frac{1+\delta}{2}\right)^{k/2}\left(\frac{2}{5}c_1-\frac{1-\delta}{25}
\left(\frac{1+\delta}{2}\right)^{k/2}\right)\nonumber\\
&\ge &(1-\delta)\left(\frac{1+\delta}{2}\right)^{k/2}\left(4-\frac{1}{25}\left(\frac{1+\delta}{2}\right)^{k/2}\right)\, ,
\end{eqnarray}
where in the third line above we used \eqref{e:epsilon-zeta} and neglected the first term in the square brackets
in the line above, in the fourth line we used that $\rho\le 1$ and again \eqref{e:epsilon-zeta}, and in the last
line we used that $c_1\ge 10$.

By combining \eqref{e:upper-A} and \eqref{e:lower-B-A} we get
\begin{eqnarray*}
B-A-3\epsilon^{k+1}A&\ge & B-A-3\epsilon^{k+1}\cdot 9c_1^2\\
&\ge &(1-\delta)\left(\frac{1+\delta}{2}\right)^{k/2}\left(4-\frac{1}{25}\left(\frac{1+\delta}{2}\right)^{k/2}\right) -27\left(\frac{1-\delta}{20c_1}\right)^{k-1}\left(\frac{1-\delta}{20c_1}\right)^2c_1^2\\
&= & (1-\delta)\left(\frac{1+\delta}{2}\right)^{k/2}\left(4-\frac{1}{25}\left(\frac{1+\delta}{2}\right)^{k/2}\right) - \frac{27}{20^2}\left(\frac{1-\delta}{20c_1}\right)^{k-1}(1-\delta)^2\\
&\ge&(1-\delta)\left(\frac{1+\delta}{2}\right)^{k/2}\left(4-\frac{1}{25}\left(\frac{1+\delta}{2}\right)^{k/2}\right)  -27\left(\frac{1-\delta}{20}\right)^{k+1}\\
&\ge & 3 (1-\delta)\left(\frac{1+\delta}{2}\right)^{k/2}-27\left(\frac{1-\delta}{20}\right)^{k+1}\\
&\ge &3 \frac12 \left(\frac12\right)^{k/2}-27\left(\frac{1}{20}\right)^{k+1}\\
&=& 3\left(\frac12\right)^{k/2+1}-27\left(\frac{1}{20}\right)^{k+1} >0
\end{eqnarray*}
for all $k\ge 1$.
In the penultimate line we used that $1-\delta\ge 1/2$, $(1+\delta)/2\ge 1/2$ and $1-\delta<1$.

The proof is now complete.
\qed

An alternative approach to the oscillation reduction in case of rotationally
stable processes is given in \cite[Lemma 8]{BKK}. We note that for non-local
operators the oscillation reduction technique seems much harder than for the
Laplacian, because the subtraction used in this process may destroy global nonnegativity.
In this regard, we note that there
is a gap of this nature in the
proof of \cite[Lemma 3.3]{CS98}.

In the remainder of this section, we assume that $Q\in \partial D$.
Since $D$ is Greenian, the Green function $G_D(x,y)$, $x,y\in D$, is well defined. Fix $x_0\in D$ and set
$$
M_D(x,y)=\frac{G_D(x,y)}{G_D(x_0,y)}\, ,\qquad x,y\in D\,, y\neq x_0.
$$
Let $r<\frac12 \min\{\mathrm{dist}(x,Q),\mathrm{dist}(x_0,Q)\}$.
Since $G_D(x,\cdot)$ and $G_D(x_0, \cdot)$ are regular harmonic in $D\cap B(Q,2r)$
and vanish in $D^c\cap B(Q,2r)$,
using Lemma \ref{l:oscillation-reduction},
one immediately gets the following.

\begin{corollary}\label{c:oscillation-reduction}
The limit $M_D(x, Q):=\lim_{D\ni y\to Q}M_D(x, y)$ exists.
Furthermore, there exist positive constants $c$ and $\beta$ depending
on $(\Psi, \gamma_1, \gamma_2, \kappa)$ such that
for any $r\in (0, \frac12(R\wedge\mathrm{dist}(x_0,Q))$, any
$y\in D\cap B(Q, r)$
and any $x\in D\setminus B(Q, 2r)$,
\begin{equation}\label{e:martin-estimate-boundary}
|M_D(x, y)-M_D(x, Q)|\le cM_D(x, A_r(Q))\left(\frac{|y-Q|}{r}\right)^\beta.
\end{equation}
\end{corollary}

\pf Put
$$
u(y):=\frac{G_D(x, y)}{G_D(x, A_r(Q))}, \qquad
v(y):=\frac{G_D(x_0, y)}{G_D(x_0, A_r(Q))}.
$$
Since $u$ and $v$ satisfy the assumptions of Lemma \ref{l:oscillation-reduction}, there exists
$$
g:=\lim_{D\ni y\to Q}\frac{u(y)}{v(y)}=\lim_{D\ni y\to Q} \frac{M_D(x,y)}{M_D(x, A_r(Q))}\, .
$$
This implies the existence of $M_D(x,Q):=\lim_{D\ni y\to Q}M_D(x,y)$. Furthermore, by \eqref{e:zoran5.6} we have that
$$
\left|\frac{u(y)}{v(y)}-g\right| \le c\left(\frac{|y-Q|}{r}\right)^{\beta}\, , \quad \text{for all }y\in D\cap B(Q,r)\, .
$$
Equivalently,
$$
\left|\frac{M_D(x,y)}{M_D(x, A_r(Q))}-\frac{M_D(x,Q)}{M_D(x, A_r(Q))}\right|\le c\left(\frac{|y-Q|}{r}\right)^{\beta}\, , \quad \text{for all }y\in D\cap B(Q,r)\, ,
$$
which is \eqref{e:martin-estimate-boundary}.
\qed

Recall that $X^D$ is the process $X$ killed upon exiting $D$.
As the process $X^D$ satisfies Hypothesis (B) in \cite{KW}, $D$  has
a Martin boundary $\partial_M D$ with respect to $X^D$ satisfying the following properties:
\begin{description}
\item{(M1)} $D\cup \partial_M D$ is
a compact metric space (with the metric denoted by $d$);
\item{(M2)} $D$ is open and dense in $D\cup \partial_M D$,  and its relative topology coincides with its original topology;
\item{(M3)}  $M_D(x ,\, \cdot\,)$ can be uniquely extended  to $\partial_M D$ in such a way that
\begin{description}
\item{(a)}
$ M_D(x, y) $ converges to $M_D(x, w)$ as $y\to w \in \partial_M D$ in the Martin topology;
\item{(b)} for each $ w \in D\cup \partial_M D$ the function $x \to M_D(x, w)$  is excessive with respect to $X^D$;
\item{(c)} the function $(x,w) \to M_D(x, w)$ is jointly continuous on $D\times (D\cup \partial_M D)$ in the Martin topology and
\item{(d)} $M_D(\cdot,w_1)\not=M_D(\cdot, w_2)$ if $w_1 \not= w_2$ and $w_1, w_2 \in \partial_M D$.
\end{description}
\end{description}

We will say that a point $w\in \partial_M D$ is a finite Martin boundary point if there exists a bounded sequence
$(y_n)_{n\ge 1}\subset D$ converging to $w$ in the Martin topology. The finite part of the Martin boundary will be denoted by $\partial_M^f D$.
Recall that a point $w$ on the Martin boundary $\partial_MD$ of $D$ is said to be associated
with $Q\in \partial D$ if there is a sequence $(y_n)_{n\ge 1}\subset D$ converging to $w$
in the Martin topology and to $Q$ in the Euclidean topology. The set of Martin
boundary points associated with $Q$ is denoted by $\partial_M^QD$.

\begin{prop}\label{p:infinite-mb}
$\partial_M^QD$ consists of exactly one point.
\end{prop}

\pf
We first note that $\partial_M^QD$ is not empty.
Indeed, let $(y_n)_{n\ge 1}\subset D$ converge to $Q$ in the Euclidean topology. Since $D\cup \partial_M D$ is a compact metric space with the Martin metric,  there exist a subsequence $(y_{n_k})_{k\ge 1}$ and $w\in D\cup \partial_M D$ such that $\lim_{k\to \infty}d(y_{n_k},w)=0$. Clearly, $w\notin D$ (since relative topologies on $D$ are equivalent). Thus we have found a sequence $(y_{n_k})_{k\ge 1}\subset D$ which converges to $w\in \partial_M D$ in the Martin topology and to $Q$ in the Euclidean topology.

Let $w\in \partial_M^QD$ and let $M_D(\cdot, w)$ be the corresponding Martin kernel.
If  $(y_n)_{n\ge 1}\subset D$ is a sequence converging to $w$
in the Martin topology and to $Q$ in the Euclidean topology,
then, by (M3)(a), $M_D(x,y_n)$ converge to $M_D(x,w)$. On the other hand, $|y_n-Q|\to 0$, thus by Corollary \ref{c:oscillation-reduction},
$$
\lim_{n\to \infty}M_D(x,y_n)=M_D(x,Q).
$$
Hence, for each $w\in \partial_M^Q D$ it holds that $M_D(\cdot, w)=M_D(\cdot, Q)$.
Since, by (M3)(d),  for two different Martin boundary points $w^{(1)}$ and $w^{(2)}$ it always holds that $M_D(\cdot, w^{(1)})\neq M_D(\cdot, w^{(2)})$, we conclude that $\partial_M^QD$ consists of exactly one point.
\qed

Because of the proposition above, we will also use $Q$ to denote the point on the Martin boundary $\partial_M^QD$
associated with $Q$.
Note that it follows from the proof of Proposition \ref{p:infinite-mb} that if $(y_n)_{n\ge 1}$ converges to $Q$ in the Euclidean topology, then it also converges to $Q$ in the Martin topology.

For $\epsilon >0$ let
\begin{equation}\label{e:definition-U_K}
K_{\epsilon}:=\left\{w\in \partial_M D: d(w, Q) \ge \epsilon\right\}
\end{equation}
be a closed subset of $\partial_M D$. By the definition of Martin boundary,
for each $w\in K_{\epsilon}$ there exists a sequence $(y_n^w)_{n\ge 1}\subset D$ such that
$\lim_{n\to \infty} d(y_n^w, w)=0$. Without loss of generality we may assume that
$d(y_n^w, w)<\frac{\epsilon}{2}$ for all $n\ge 1$.

\begin{lemma}\label{l:boundedness-U-K}
There exists $c=c(\epsilon)>0$ such that $|y_n^w-Q|\ge c$
for all $w\in K_{\epsilon}$ and all $n\ge 1$.
\end{lemma}
\pf
Suppose the lemma is not true. Then $\{y_n^w:\, w\in K_{\epsilon}, n\in \N\}$
contains a sequence $(y_{n_k}^{w_k})_{k\ge 1}$ such that $\lim_{k\to \infty}|y_{n_k}^{w_k}-Q|= 0$.
Then also
$\lim_{k\to \infty}d(y_{n_k}^{w_k}, Q)= 0$. On the other hand,
$$
d(y_{n_k}^{w_k}, Q)\ge d(w_k, Q)-d(y_{n_k}^{w_k}, w_k)\ge \epsilon-\frac{\epsilon}{2}=\frac{\epsilon}{2}\, .
$$
This contradiction proves the claim. \qed

We continue by showing that $M_D(\cdot, Q)$ is harmonic in $D$ with respect to $X$.

\begin{lemma}\label{l:mk-integrability}
For every
bounded open $U\subset \overline{U}\subset D$ and every $x\in D$, $M_D(X_{\tau_U}, Q)$ is $\P_x$-integrable.
\end{lemma}
\pf Let $(y_m)_{m\ge 1}$ be a sequence in $D\setminus \overline{U}$ such that $|y_m-Q|\to 0$. Then $M_D(\cdot, y_m)$ is regular harmonic in $U$. Hence, by Fatou's lemma and Corollary \ref{c:oscillation-reduction},
\begin{eqnarray*}
  &&  \E_x[M_D(X_{\tau_U}, Q)]\,=\, \E_x[\lim_{m\to \infty}M_D(X_{\tau_U}, y_m)]
    \,\le \,\liminf_{m\to \infty}\E_x[M_D(X_{\tau_U}, y_m)]\\
    &&=\,\liminf_{m\to \infty} M_D(x,y_m)\,=\,M_D(x, Q)<\infty\, .
\end{eqnarray*}
\qed

Using the results above, we can get the following result.

\begin{lemma}\label{l:mk-harmonic}
Suppose that $D$ is either (i) a bounded open set or (ii) an unbounded  open set
and $X$ is transient.
For any $x\in D$ and $\rho\in (0, R\wedge(\delta_D(x)/2)]$,
$$
    M_D(x, Q)=\E_x[M_D(X_{\tau_{B(x,\rho)}},
Q)]\, .
$$
\end{lemma}

\pf Fix $x\in D$ and a positive $r<R\wedge \frac{\delta_D(x)}2$. Let
$$
\eta_m:=\left(\frac{\kappa}2\right)^mr \quad \mbox{and  }\ z_m=A_{\eta_m}(Q),
\quad m=0, 1, \dots.
$$
Note that
$$
B(z_m, \eta_{m+1})\subset D\cap B(Q, \frac{\eta_m}2)\subset D\cap B(Q,\eta_m)\subset
D\cap B(Q, r)\subset D\setminus B(x, r)
$$
for all $m\ge 0$. Thus by the harmonicity of $M_D(\cdot, z_m)$, we have
$$
M_D(x, z_m)=\E_x\left[M_D(X_{\tau_{B(x, r)}}, z_m)\right].
$$

On the other hand, by Theorem \ref{ubhp}, there exist $m_0=m_0(\kappa)\ge 2$
and $c_1=c_1(\Psi, \gamma_1, \gamma_2, \kappa)>0$
such that for every $w\in D\setminus B(Q, \eta_m)$ and $y\in D\cap B(Q, \eta_{m+1})$,
$$
M_D(w, z_m)=\frac{G_D(w, z_m)}{G_D(x_0, z_m)}\le c_1\frac{G_D(w, y)}{G_D(x_0, y)}=c_1M_D(w, y), \qquad m\ge m_0.
$$
Letting $y\to Q$ we get
\begin{equation}\label{e:KSV15.8}
M_D(w, z_m)\le c_1M_D(w, Q), \qquad m\ge m_0, w\in D\setminus B(Q, \eta_m).
\end{equation}

To prove this lemma, it suffices to show that $\{M_D(X_{\tau_{B(x, r)}}, z_m): m\ge m_0\}$ is $\P_x$-uniformly
integrable. Since $M_D(X_{\tau_{B(x, r)}}, Q)$ is integrable by  Lemma \ref{l:mk-integrability}, for
any $\epsilon>0$, there is an $N_0>1$ such that
\begin{equation}\label{e:KSV15.9}
\E_x\left[M_D(X_{\tau_{B(x, r)}}, Q); M_D(X_{\tau_{B(x, r)}}, Q)>N_0/c_1\right]<\frac{\epsilon}{2c_1}.
\end{equation}
Note that by \eqref{e:KSV15.8} and \eqref{e:KSV15.9},
\begin{eqnarray*}
&&\E_x\left[M_D(X_{\tau_{B(x, r)}}, z_m); M_D(X_{\tau_{B(x, r)}}, z_m)>N_0 \mbox{ and }
 X_{\tau_{B(x, r)}}\in D\setminus B(Q, \eta_m)\right]\\
&&\le c_1\E_x\left[M_D(X_{\tau_{B(x, r)}}, Q); c_1M_D(X_{\tau_{B(x, r)}}, Q)>N_0\right]
<c_1\frac{\epsilon}{2c_1}=\frac{\epsilon}{2}.
\end{eqnarray*}
By \eqref{e:KSV7-4.7}, we have for $m\ge m_0$,
\begin{eqnarray*}
&&\E_x\left[M_D(X_{\tau_{B(x, r)}}, z_m);X_{\tau_{B(x, r)}}\in D\cap B(Q, \eta_m)\right]\\
&& =\int_{D\cap B(Q, \eta_m)}M_D(w, z_m)K_{B(x, r)}(x, w)dw\\
&&\le c_2 \phi(r^{-2})^{-1}\int_{D\cap B(Q, \eta_m)}M_D(w, z_m)j(|w-x|-r)dw
\end{eqnarray*}
for some $c_2=c_2(\Psi, \gamma_1, \gamma_2)>0$.
Since $|w-x|\ge |x-Q|-|Q-w|\ge \delta_D(x)-\eta_m\ge
\frac{7}4r$, using \eqref{e:doubling-condition} and \eqref{e:asmpbofjat0}, we get that
\begin{eqnarray}
&&\E_x\left[M_D(X_{\tau_{B(x, r)}}, z_m);X_{\tau_{B(x, r)}}\in D\cap B(Q, \eta_m)\right]\nonumber\\
&&\le c_3 j(r)\phi(r^{-2})^{-1}\int_{D\cap B(Q, \eta_m)}M_D(w, z_m)dw\nonumber\\
&&\le c_4r^{-d}\int_{D\cap B(Q, \eta_m)}M_D(w, z_m)dw\nonumber\\
&&=c_4r^{-d}G_D(x_0, z_m)^{-1}\int_{D\cap B(Q, \eta_m)}G_D(w, z_m)dw
\label{e:KSV15.10}
\end{eqnarray}
for some $c_3=c_3(\Psi, \gamma_1, \gamma_2)>0$ and $c_4=c_4(\Psi, \gamma_1, \gamma_2)>0$.
Note that, by Lemma \ref{l:lemma3.1},
\begin{equation}\label{e:KSV15.11}
G_D(x_0, z_m)^{-1}\le C_3(\phi(\eta_{m_0}^{-2}))^{-1}
  \xi^{m-m_0}
  \phi(\eta_m^{-2})
G_D(x_0, z_{m_0})^{-1}.
\end{equation}
By \eqref{e:H1}, there exists $c_5=c_5(\Psi, \gamma_1, \gamma_2)>0$
such that for any $\eta<1$,
$$
\int^{\eta}_0\frac1{s\phi(s^{-2})}ds\le c_5\phi(\eta^{-2})^{-1}.
$$
Thus by Lemma \ref{dec} in case $D$ is bounded and by Theorem \ref{t:greenh} in case of unbounded $D$,
\begin{equation}\label{e:KSV15.12}
\int_{B(Q, \eta_m)}G_D(w, z_m)dw\le c_6\int_{B(z_m, 2\eta_m)}\frac1{|w-z_m|^d\phi(|w-z_m|^{-2})}dw
\le c_7 \phi((2\eta_m)^{-2})^{-1}
\end{equation}
for some constants $c_6=c_6(\Psi, \gamma_1, \gamma_2)>0$ and $c_7=c_7(\Psi, \gamma_1, \gamma_2)>0$.
It follows from \eqref{e:KSV15.10}--\eqref{e:KSV15.12} that
\begin{eqnarray*}
&&\E_x\left[M_D(X_{\tau_{B(x, r)}}, z_m);X_{\tau_{B(x, r)}}\in D\cap B(Q, \eta_m)\right]\\
&&\le c_8r^{-d}(\phi(\eta_{m_0}^{-2}))^{-1}G_D(x_0, z_{m_0})^{-1}\frac{\phi(\eta_m^{-2})}
{\phi((2\eta_m)^{-2})}
\xi^{m-m_0}.
\end{eqnarray*}
Applying Lemma \ref{l:phi-property}, we get
\begin{eqnarray*}
\E_x\left[M_D(X_{\tau_{B(x, r)}}, z_m);X_{\tau_{B(x, r)}}\in D\cap B(Q, \eta_m)\right]
\le c_9r^{-d}(\phi(\eta_{m_0}^{-2}))^{-1}G_D(x_0, z_{m_0})^{-1}
  \xi^{m-m_0}.
\end{eqnarray*}
Thus there exists $N>0$ such that for all $m\ge N$,
$$
\E_x\left[M_D(X_{\tau_{B(x, r)}}, z_m);X_{\tau_{B(x, r)}}\in D\cap B(Q, \eta_m)\right]\le \frac{\epsilon}2.
$$
Consequently, for all $m\ge N$,
$$
\E_x\left[M_D(X_{\tau_{B(x, r)}}, z_m);M_D(X_{\tau_{B(x, r)}}, z_m)>N\right]\le \epsilon,
$$
which implies that $\{M_D(X_{\tau_{B(x, r)}}, z_m): m\ge m_0\}$ is $\P_x$-uniformly
integrable.
\qed

Using this, we can easily get the following

\begin{thm}\label{t:mk-harmonic}
The function $M_D(\cdot, Q)$ is harmonic in $D$ with respect to $X$.
\end{thm}

\pf The proof is basically the same as that of \cite[Theorem 3.9]{KSV6}. We write the details here for
completeness.  Let $h(x):= M_D(x, Q)$. Consider a relatively compact open set $D_1\subset\overline{D_1}
\subset D$, and put $r(x)=R\wedge (\frac13\delta_D(x))$ and $B(x)=B(x, r(x))$.
Define a sequence $\{T_m: m\ge 1\}$ of stopping times as follows: $T_1:=\inf\{t>0: X_t\notin B(X_0)\}$,
and for $m\ge 2$,
$$
T_m:=\begin{cases} T_{m-1}+\tau_{B(X_{T_{m-1}})}\cdot \theta_{T_{m-1}} & \mbox{ if } X_{T_{m-1}}\in D_1\\
\tau_{D_1} & \mbox{ otherwise}.
\end{cases}
$$
Note that $X_{\tau_{D_1}}\in \partial D_1$ on $\cap^\infty_{n=1}\{T_n<\tau_{D_1}\}$. Thus, since
$\lim_{m\to\infty}T_m=\tau_{D_1}$ $\P_x$-a.s. and $h$ is continuous in $D$, using the quasi-left continuity
of $X^D$, we have $\lim_{m\to\infty}h(X^D_{T_m})=h(X^D_{\tau_{D_1}})$ on $\cap^\infty_{n=1}\{T_n<\tau_{D_1}\}$.
Now by the dominated convergence theorem and Lemma \ref{l:mk-harmonic},
\begin{eqnarray*}
h(x)&=&\lim_{m\to\infty}\E_x[h(X^D_{T_m}); \cup^\infty_{n=1}\{T_n=\tau_{D_1}\}]
+\lim_{m\to\infty}\E_x[h(X^D_{T_m}); \cap^\infty_{n=1}\{T_n<\tau_{D_1}\}]\\
&=&\E_x[h(X^D_{\tau_{D_1}}); \cup^\infty_{n=1}\{T_n=\tau_{D_1}\}]+
\E_x[h(X^D_{\tau_{D_1}}); \cap^\infty_{n=1}\{T_n<\tau_{D_1}\}]\\
&=&\E_x[h(X^D_{\tau_{D_1}})].
\end{eqnarray*}
\qed

Part (b) of the following result is proved in \cite[Lemma 4.18]{KSV9}. Part (a) is even simpler.

\begin{lemma}\label{l:irregular}
(a) Let $D$ be a bounded open set and
suppose that $u$ is a bounded nonnegative harmonic function for $X^D$. If there exists a polar set $N\subset \partial D$ such that for any $z\in \partial D\setminus N$
\begin{equation}\label{e:not-polar}
\lim_{D\ni x\to z} u(x)=0\, ,
\end{equation}
then $u$ is identically equal to zero.

\noindent
(b) Let $D$ be an unbounded open set and
suppose that $u$ is a bounded nonnegative harmonic function for $X^D$. If there exists a polar set $N\subset \partial D$ such that for any $z\in \partial D\setminus N$ \eqref{e:not-polar} holds true and additionally
$$
\lim_{D\ni x\to \infty} u(x)=0\, ,
$$
then $u$ is identically equal to zero.
\end{lemma}

The next result completes the proof of Theorem \ref{t:main-theorem}.
Recall that a point $z\in \partial D$ is said to regular boundary point of $D$ if $\P_z(\tau_D=0)=1$
and an irregular boundary point if $\P_z(\tau_D=0)=0$. The set of irregular boundary points is polar.

\begin{thm}\label{t:m}
Assume that either $D$ is bounded, or $D$ is unbounded and $X$ is transient.
Then $Q$ is a minimal Martin boundary point, that is,
$M_D(\cdot, Q)$ is a minimal harmonic function.
\end{thm}

\pf Let $h$ be a positive harmonic function for  $X^{D}$
such that $h\le M_D(\cdot, Q)$. By
the Martin representation in \cite{KW},
there is a finite measure on $\partial_M D$ such that
$$
    h(x)=\int_{\partial_M D}M_D(x,w)\, \mu(dw)=\int_{\partial_M D\setminus\{Q\}}M_D(x,w)\, \mu(dw)+ M_D(x,
Q)\mu(\{Q\})\, .
$$
In particular, $\mu(\partial_M D)=h(x_0)\le
M_D(x_0, Q)=1$ (because of the normalization at $x_0$). Hence, $\mu$ is a sub-probability measure.

For $\epsilon >0$, $K_{\epsilon}$
is the compact subset of $\partial_M D$ defined in \eqref{e:definition-U_K}.
Define
\begin{equation}\label{d:definition-u}
    u(x):=\int_{ K_{\epsilon} }M_D(x,w)\, \mu(dw).
\end{equation}
Then $u$ is a positive harmonic function with respect to  $X^{D}$
satisfying
\begin{align}\label{e:newm1}
u(x)\le h(x)-\mu(\{Q\})M_D(x,
Q)\le \big(1-\mu(\{Q\})\big)M_D(x, Q)\, .
\end{align}

Let $c=c(\epsilon)>0$ be the constant from Lemma \ref{l:boundedness-U-K}. Hence, for $w\in K_{\epsilon}$ and $(y_n^w)_{n\ge 1}$ a sequence such that $\lim_{n\to \infty}d(y_n^w,w)=0$, it holds that $|y_n^w-Q|\ge c$. Fix $x_1\in D\cap B(Q,c/2)$ and choose arbitrary $y_0\in D\setminus B(Q,c)$. For any $x\in D\cap B(Q,c/2)$ and any $y\in D\setminus B(Q,c)$ we have that
$$
\frac{G_D(x,y)}{G_D(x_0,y)}=\frac{G_D(x,y)}{G_D(x_1,y)}\, \frac{G_D(x_1,y)}{G_D(x_0,y)}\le c_1 \frac{G_D(x,y_0)}{G_D(x_1,y_0)}\, \frac{G_D(x_1,y)}{G_D(x_0,y)}\, .
$$
Here the inequality follows from Theorem \ref{ubhp} applied to functions $G_D(\cdot, y)$ and $G_D(\cdot, y_0)$ which are regular harmonic in $D\cap B(Q,c)$ and vanish in $D\setminus B(Q,c)$. Now fix $w\in K_{\epsilon}$ and apply the above inequality to $y_n^w$ to get
\begin{eqnarray*}
M_D(x,w)&=&\lim_{n\to \infty}\frac{G_D(x,y_n^w)}{G_D(x_0, y_n^w)}\le c_1 \frac{G_D(x,y_0)}{G_D(x_1,y_0)}\, \lim_{n\to
\infty}\frac{G_D(x_1, y_n^w)}{G_D(x_0, y_n^w)}\\
&=&c_1  \frac{G_D(x,y_0)}{G_D(x_1,y_0)}\, M_D(x_1,w)\le c_1  \frac{G_D(x,y_0)}{G_D(x_1,y_0)} \sup_{w\in K_{\epsilon}}M_D(x_1,w)\\
&\le &c_2  \frac{G_D(x,y_0)}{G_D(x_1,y_0)} = c_3  G_D(x,y_0)\, .
\end{eqnarray*}
In the  last line we used property (M3) (c) of the Martin kernel. Thus,
\begin{equation}\label{e:bound-on-K-epsilon}
M_D(x,w)\le c_3 G_D(x,y_0)\, , \qquad x\in D\cap B(Q,c/2), w\in K_{\epsilon}\, .
\end{equation}

Choose $r<c/4$. For any $x\in D\setminus B(Q, 2r)$ and $y\in D\cap B(Q, r)$ with $r$ small enough, by Theorem \ref{ubhp} applied to $G_D(x, \cdot)$ and $G_D(x_0,\cdot)$, we have
$$
\frac{G_D(x, y)}{G_D(x_0, y)}\le c_4\frac{G_D(x, A_r(Q))}{G_D(x_0, A_r(Q))}.
$$
Letting $D\ni y\to Q$, we get
\begin{equation}\label{e:martinkernelub}
M_D(x, Q)\le c_4\frac{G_D(x, A_r(Q))}{G_D(x_0, A_r(Q))}=c_5 G_D(x, A_r(Q))\, ,\quad x\in  D\setminus B(Q, 2r)\, .
\end{equation}
Recall that $\lim_{D\ni x\to z}G_D(x,y)=0$ for every regular $z\in \partial D$.
Since $r<c/4$ can be arbitrarily small, we see from \eqref{e:martinkernelub} and
 \eqref{e:newm1} that  $\lim_{D\ni x, x\to z}u(x)=0$ for every regular $z\in \partial D$, $z\neq Q$.

Assume $D$ is bounded. Fix $r<c/4$. It follows from Lemma \ref{dec} that for all $x\in D\setminus B(Q,2r)$,
\begin{equation}\label{e:martinkernelub3}
G_D(x, A_r(Q))\le c_6 \frac{\Phi(|x-A_r(Q)|)}{|x-A_r(Q)|^d}\le c_6\sup_{a\ge r} \frac{\Phi(a)}{a^d}\le c_7\, .
\end{equation}
From \eqref{e:martinkernelub} and \eqref{e:newm1}  we conclude that $u$ is
bounded in $x\in D\setminus B(Q,2r)$. Similarly, for every $x\in D\cap B(Q,c/2)$
we have that  $G_D(x,y_0)\le c_8 \sup_{a\ge c/2}\Phi(c)c^{-d}=:c_9$
(recall $y_0\in D\setminus B(Q,c)$). Hence by \eqref{e:bound-on-K-epsilon}
and  \eqref{e:newm1} we see that $u$ is bounded on $D\cap B(Q,c/2)$.
Thus $u$ is bounded on $D$. Now it follows from Lemma \ref{l:irregular} (a) that $u\equiv 0$ in $D$.

If $D$ is unbounded, we argue as follows. Since $G_D(x,A_r(Q))\le G(x,A_r(Q))$,
it follows from \eqref{e:bound-on-K-epsilon} and Lemma \ref{l:gfnatinfty} that
$\lim_{D\ni x\to \infty} M_D(x,Q)=0$. Hence by \eqref{e:newm1}
$\lim_{D\ni x\to \infty}u(x)=0$. Thus, there exists $R\ge 2$ such that $u(x)\le 1$
for all $x\in D\setminus B(Q,R)$. Fix $r<c/4\wedge 1$ and let $x\in D\cap (B(Q,R)\setminus B(Q,2r))$. By \eqref{e:martinkernelub} and Theorem \ref{t:greenh},
$$
M_D(x,Q)\le c_5 G(x, A_r(Q))\le c_5 C_2(R) \frac{\Phi(|x-A_r(Q)|)}{|x-A_r(Q)|^d}\le c_{10}\sup_{a\ge r} \frac{\Phi(a)}{a^d}\le c_{11}\, .
$$
It follows that $u$ is bounded in $D\cap (B(Q,R)\setminus B(Q,2r))$. The proof that $u$ is bounded on $D\cap B(Q,c/2)$ is the same as in the case of a bounded $D$. Hence, $u$ is bounded, and again we conclude from Lemma \ref{l:irregular} (b) that $u\equiv 0$ in $D$.

We see from \eqref{d:definition-u} that $\nu=\mu_{| K_{\epsilon}}=0$. Since
$\epsilon >0$ was arbitrary and $\partial_M D\setminus\{Q\}=\cup_{\epsilon >0} K_{\epsilon}$,
we see that $\mu_{|\partial_M D\setminus\{Q\}}=0$. Hence $h=\mu(\{Q\})
M_D(\cdot, Q)$ showing that $
M_D(\cdot, Q)$ is minimal.
\qed

In the next two results we assume that $D$ is a $\kappa$-fat set.
Then one can define $\Xi:\partial D \to \partial_M^f D$
so that $\Xi(Q)$ is the unique element of $\partial_M^QD$, cf.~Proposition \ref{p:infinite-mb}.

\begin{thm}\label{mbid}
Suppose that either $D$ is bounded, or $D$ is unbounded and $X$ in transient. If $D$ is a
$\kappa$-fat set, then the finite part of the minimal Martin boundary of $D$ and the finite
part of the  Martin boundary of $D$ both coincide with  the Euclidean boundary $\partial D$ of $D$.
More precisely, $\Xi$ is 1-1- and onto.
\end{thm}
\pf
Since every finite Martin boundary point is associated with some $Q\in \partial D$, we see that $\Xi$ is onto.
We show now that $\Xi$ is 1-1.
If not, there are $Q, Q'\in \partial D$, $Q\neq Q'$, such that $\Xi(Q)=\Xi(Q')=w$.
Then $M_D(\cdot, Q)=M_D(\cdot, w)= M_D(\cdot, Q')$. Choose $r>0$ small enough and
satisfying  $r<|Q-Q'|/4$. By \eqref{e:martinkernelub} and \eqref{e:martinkernelub3}
we see that there exists a constant $c_1=c_1(Q)$ such that $M_D(x, Q)\le c_1$ for
all $x\in D\setminus B(Q,2r)$. Similarly, there exists $c_2=c_2(Q')$  such that
$M_D(x, Q')\le c_2$ for all $x\in D\setminus B(Q',2r)$. Since $B(Q, 2r)$ and
$B(Q',2r)$ are disjoint, we conclude that $M_D(\cdot, Q)=M_D(\cdot, Q')$ is
bounded  on $D$ by $c_1\vee c_2$. Again by \eqref{e:martinkernelub},
$\lim_{D\ni x\to z}M_D(x,Q)=0$ for all regular $z\in \partial D$.
In case of unbounded $D$, we showed in the proof of Theorem \ref{t:m}
that $\lim_{x\to \infty}M_D(x,Q)=0$. Hence by Lemma \ref{l:irregular}
we conclude that $M_D(\cdot, Q)\equiv 0$. This is a contradiction with $M_D(x_0, Q)=1$.

The statement about the minimal Martin boundary follows from Theorem \ref{t:m}.
\qed

As a consequence of the result above and the general result of \cite{KW}), we have the following
Martin representation for nonnegative harmonic functions with respect to the killed
process $X^D$.

\begin{thm}
Suppose that $D$ is a bounded $\kappa$-fat set.
Then $\Xi:\partial D \to \partial_M D$ is a homeomorphism. Furthermore,
for any nonnegative function $u$
which is harmonic with respect to $X^D$, there exists a unique finite measure
$\mu$ on $\partial D$ such that
$$
u(x)=\int_{\partial D}M_D(x, z)\mu(dz), \qquad x\in D.
$$
\end{thm}
\pf
Let $Q\in \partial D$ and $x\in D$. Choose $r<\frac12 \min\{R,\mathrm{dist}(x,Q),\mathrm{dist}(x_0,Q)\}$ so that $x\in D\setminus B(Q, 2r)$. Let $Q'\in \partial D\cap B(Q,r/2)$. Since $D$ is $\kappa$-fat at $Q'$, by Corollary \ref{c:oscillation-reduction} there exists $M_D(x,Q')=\lim_{D\ni y\to Q'}M_D(x, y)$. Further, by letting $y\to Q'$ in \eqref{e:martin-estimate-boundary} we get that
$$
|M_D(x, Q')-M_D(x, Q)|\le cM_D(x, A_r(Q))\left(\frac{|Q'-Q|}{r}\right)^\beta.
$$
This shows
that if $(Q_n)_{n\ge 1}$ is a sequence of points
in $\partial D$ converging to $Q\in \partial D$, then $M_D(\cdot,Q)=\lim_{n\to \infty}M_D(\cdot, Q_n)$.

In order to show that $\Xi$ is continuous we proceed as follows.
Let $Q_n\to Q$ in $\partial D$. Since $\partial_M D$ is compact, $(\Xi(Q_n))_{n\ge 1}$ has a subsequence $(\Xi(Q_{n_k}))_{k\ge 1}$ converging in the Martin topology to some $w\in \partial_M D$. By property (M3), $M_D(\cdot, \Xi(Q_{n_k}))\to M_D(\cdot, w)$. On the
other hand, by the first part of the proof,
$M_D(\cdot, \Xi(Q_{n_k}))=M_D(\cdot, Q_{n_k})\to M_D(\cdot, Q)$, implying that $w=\Xi(Q)$. This shows in fact that $(\Xi(Q_n))_{n\ge 1}$ is convergent with the limit $\Xi(Q)$. Using the fact that $\partial D$ is compact, the proof of the continuity of the inverse is similar.

The Martin representation for nonnegative harmonic functions is now a consequence of the general result form \cite{KW}.
\qed

\bigskip

{\bf Acknowledgements.} We thank the referee for helpful comments on the
first version of this paper.

\vspace{.1in}
\begin{singlespace}

\end{singlespace}

\end{doublespace}
\vskip 0.1truein

{\bf Panki Kim}

Department of Mathematical Sciences and Research Institute of Mathematics,

Seoul National University, Building 27, 1 Gwanak-ro, Gwanak-gu Seoul 151-747, Republic of Korea

E-mail: \texttt{pkim@snu.ac.kr}

\bigskip

{\bf Renming Song}

Department of Mathematics, University of Illinois, Urbana, IL 61801,
USA

E-mail: \texttt{rsong@math.uiuc.edu}

\bigskip

{\bf Zoran Vondra\v{c}ek}

Department of Mathematics, University of Zagreb, Zagreb, Croatia

Email: \texttt{vondra@math.hr}

\small
\begin{thebibliography}{99}

\bibitem{BF} C. Berg and G. Forst: {\em Potential Theory on
Locally Compact Abelian Groups}. Springer, 1975.

\bibitem{B} K. Bogdan:  The boundary Harnack principle for the
fractional Laplacian.
 {\em Studia Math.}  {\bf 123(1)}(1997), 43--80.

\bibitem{B99} K. Bogdan: Representation of $\alpha$-harmonic functions in Lipschitz domains.
{\em Hiroshima Math. J.} {\bf 29} (1999), 227--243.

\bibitem{BGR1} K. Bogdan, T. Grzywny and M. Ryznar:
Density and tails of unimodal convolution semigroups.
{\em J. Funct. Anal.} {\bf 266} (2014), 3543--3571.

\bibitem{BGR2} K. Bogdan, T. Grzywny and M. Ryznar:
Barriers, exit time and survival probability for unimodal L\'evy processes. Preprint.
arXiv:1307.0270 [math.PR]

\bibitem{BKK} K. Bogdan, T. Kulczycki, M. Kwa\'{s}nicki:
Estimates and structure of $\alpha$-harmonic functions. \emph{Probab. Theory Rel. Fields.} {\bf 140} (2008), 345--381.

\bibitem{CKK2}  Z.-Q. Chen, P. Kim and  T. Kumagai:
On heat kernel estimates and parabolic Harnack inequality for jump
processes on metric measure spaces.   {\it Acta Mathematica Sinica,
English Series \bf 25} (2009), 1067--1086.

\bibitem{CKS2} Z.-Q. Chen, P. Kim and R. Song:
Dirichlet heat kernel estimates for fractional Laplacian with gradient perturbation.
\emph{Ann. Probab.} {\bf 40} (2012), 2483--2538.

\bibitem{CKS} Z.-Q. Chen, P. Kim and R. Song: Dirichlet heat kernel estimates for rotationally
   symmetric L\'evy processes. To appear in  {\it Proc. London Math. Soc.}

\bibitem{CK} Z.-Q. Chen and  T. Kumagai:
Heat kernel estimates for jump processes of mixed types on metric measure spaces. {\it Probab. Theory Relat. Fields}, {\bf 140} (2008), 277--317.

\bibitem{CS98} Z.-Q. Chen and R. Song:
Martin boundary and integral representation for harmonic functions of symmetric stable processes.
{\em J. Funct. Anal.} {\bf 159} (1998), 267--294.


\bibitem{Chu} K.-L. Chung: {\it Lectures from Markov Processes to
Brownian Motion}. Springer, New York (1982)



\bibitem{G} T. Grzywny:
On Harnack inequality and H\"older regularity for isotropic unimodal L\'evy processes.
To appear in \emph{Potential Anal.}

\bibitem{HW} R. A. Hunt and R. L. Wheeden: Positive
harmonic  functions on Lipschitz domains.
{\em Trans. Amer. Math. Soc.} {\bf 147} (1970), 507--527.

\bibitem{KS} P. Kim and R. Song: Two-sided estimates on the
density of Brownian motion with singular drift. {\em Illinois J.
Math.}, {\bf 50 (3)}, (2006), 635--688.

\bibitem{KS08} P. Kim and R. Song: Boundary behavior of harmonic functions for truncated
stable processes. {\em J. Theoret. Probab.} {\bf 21} (2008), 287--321.

\bibitem{KSV1} P.~Kim, R.~Song and Z.~Vondra\v{c}ek: Boundary Harnack principle for subordinate Brownian motions. \emph{Stoch.~Processee Appl.} {\bf 119} (2009), 1601--1631.

\bibitem{KSV3} P.~Kim, R.~Song and Z.~Vondra\v{c}ek: Potential theory of subordinate Brownian motions
revisted. {\em  Analysis applications to finance,  essays in honour of Jia-an Yan}. Interdisciplinary
Mathematical Sciences - Vol. 13, World Scientific, 2012, 277--317.

\bibitem{KSV6} P.~Kim, R.~Song and Z.~Vondra\v{c}ek: Minimal thinness for subordinate Brownian motion in half space.   \emph{Ann.~Inst.~Fourier} {\bf 62 (3)} (2012), 1045--1080.

\bibitem{KSV2}
P.~Kim, R.~Song and Z.~Vondra\v{c}ek:  Two-sided {G}reen function estimates for killed subordinate
  {B}rownian motions. {\it Proc. London Math. Soc.} \textbf{104} (2012), 927--958.

\bibitem{KSV7} P.~Kim, R.~Song and Z.~Vondra\v{c}ek: Uniform boundary Harnack principle for rotationally
symmetric L\'evy processes in general open sets. \emph{Science China Math.} {\bf 55} (2012), 2317--2333.

\bibitem{KSV8} P.~Kim, R.~Song and Z.~Vondra\v{c}ek:  Global uniform boundary Harnack principle with explicit decay rate and its application. \emph{Stoch. Proc. Appl.}  {\bf 124(1)} (2014), 235--267.

\bibitem{KSV9} P.~Kim, R.~Song and Z.~Vondra\v{c}ek: Boundary Harnack principle and Martin boundary at infinity for subordinate Brownian motions. To appear in {\it Potential Analysis}

\bibitem{KMR} M. Kwa\'snicki, J. Ma{\l}ecki, M. Ryznar:  Suprema of L\'evy
processes. {\em Ann.~Probab.} {\bf 41} (2013) 2047--2065.

\bibitem{KW} H.~Kunita and T.~Watanabe: Markov processes and Martin boundaries I.
\emph{Illinois J.~Math.}
 {\bf 9(3)}  (1965) 485--526.


\bibitem{martin} R.S. Martin: Minimal harmonic functions.
{\em Trans. Amer. Math. Soc.} {\bf 49} (1941), 137--172.

\bibitem{MS00} K. Michalik and K. Samotij: Martin representation for $\alpha$-harmonic functions.
{\em Probab. Math. Statist. } {\bf 20} (2000), 75--91.

\bibitem{SV} R.~Song, Z.~Vondra\v{c}ek: \emph{Potential theory of subordinate Brownian motion}. In: Potential Analysis of Stable Processes and its Extensions, P. Graczyk, A. Stos, editors, Lecture Notes in Mathematics 1980, (2009) 87--176.

\bibitem{SW99} R. Song and J. Wu: Boundary Harnack principle for symmetric stable processes.
{\em J. Funct. Anal.} {\bf 168} (1999), 403--427.

\bibitem{Sz00} P. Sztonyk: On harmonic measure for L\'evy processes. {\em Probab. Math. Statist.}
{\bf 20} (2000), 383--390.

\bibitem{Sz} P. Sztonyk: Boundary potential theory for stable L\'evy processes.
{\em Colloq. Math.} {\bf 95} (2003), 191--206.
\end{thebibliography}
\end{document}